\documentclass[a4paper,10pt]{article}

\usepackage{geometry}
    \usepackage{amsmath}
    \usepackage{amsfonts,color,amssymb,mathrsfs,stmaryrd}
    \usepackage{amsthm}
    \usepackage{amssymb}
		\usepackage{url}

\usepackage{bera}

    \usepackage[normalem]{ulem}

\usepackage{graphicx}
\usepackage{wrapfig}
\usepackage{mathrsfs}

\numberwithin{equation}{section}

    \newtheorem{theo}{Theorem}\numberwithin{theo}{section}

    \newtheorem{prop}[theo]{Proposition}
    \newtheorem{lemm}[theo]{Lemma}

		\theoremstyle{remark}
    \newtheorem{rema}{Remark}
		
        \def\N{\mathbb{N}}
    \def\E{\mathbb{E}}
    \def\0{{\bf 0}}

    \renewcommand{\E}{\mathbb E \,}

    \newcommand{\Var}{{\rm Var}}

    \def\bdm{\begin{displaymath}}
    \newcommand{\edm}{\end{displaymath}}
    \def\benu{\begin{enumerate}}
    \def\eenu{\end{enumerate}}
    \def\beqn{\begin{equation}}
    \def\eeqn{\end{equation}}
    \def\be{\begin{equation}}
    \def\ee{\end{equation}}
    \def\bea{\begin{eqnarray}}
    \def\eea{\end{eqnarray}}
    \newcommand{\bean}{\begin{eqnarray*}}
    \newcommand{\eean}{\end{eqnarray*}}
    \newcommand{\bear}{\begin{eqnarray}}
    \newcommand{\eear}{\end{eqnarray}}

    \def\R{\mathbb{R}}

    \def\qed{\hfill\hbox{${\vcenter{\vbox{
        \hrule height 0.4pt\hbox{\vrule width 0.4pt height 6pt
        \kern5pt\vrule width 0.4pt}\hrule height 0.4pt}}}$}}

\def\qedskip{\smallskip\noindent}
\def\qed{\hfill $\Box$ \qedskip}

\usepackage{titlesec}

\titleformat*{\section}{\normalfont\large\bfseries}
\titleformat*{\subsection}{\normalfont\bfseries}

\date{\vspace{-0.95cm}}

\begin{document}

\title{Non-standard limits for a family of \\
autoregressive stochastic sequences}

\author{Sergey Foss\footnotemark[1] \ \ and \ Matthias Schulte\footnotemark[2]}


\date{\today}
\maketitle

\footnotetext[1]{Heriot-Watt University, Edinburgh, UK and Novosibirsk State University and Sobolev Institute of Mathematics, Russia. Email: 
    sergueiorfoss25@gmail.com.  
    The work is supported in part by Mathematical Center in Akademgorodok under agreement No. 075-15-2019-1675 with the Ministry of Science and Higher Education of the Russian Federation.}

\footnotetext[2]{Hamburg University of Technology, Germany and Heriot-Watt University, Edinburgh, UK. Email: 
    matthias.schulte@tuhh.de.}

\begin{abstract}
We consider a family of multivariate autoregressive stochastic sequences that restart when hit  a neighbourhood of the origin, and study their distributional limits when the autoregressive coefficient tends to one, the noise scaling parameter tends to zero, and the neighbourhood size varies. We obtain a non-standard limit theorem where the limiting distribution is a mixture of an atomic distribution and an absolutely continuous distribution whose marginals, in turn, are mixtures of distributions of signed absolute values of normal random variables. In particular, we provide conditions for the limiting distribution to be normal, like in the
case without restart mechanism. The main theorem is accompanied by a number of examples and auxiliary results of their own interest.
\end{abstract}

\vskip6pt
\noindent {\bf Keywords:} 
Autoregressive model, characteristic function, existence of moments, limiting distribution, normal distribution, regenerative cycle, restart mechanism, stationary distribution.

\vskip6pt
\noindent {\bf AMS 2020 Subject Classification:} Primary 60F05; Secondary 60E99; 60G50; 60J05.

\section{Introduction}

There is a permanent interest in autoregressive models for at least half a century,
both from theoretical and practical viewpoints, going back
to, say, Kesten \cite{Kes} and Vervaat \cite{Vervaat1979}, see also
\cite{Bra, Emb} and references therein.  We like to mention briefly various directions of recent research on  uni- and multivariate AR processes. 
There is a large number of papers in probability and econometrics on ``almost non-stationary'' AR processes, see e.g.\ the overview paper \cite{Rah} and the earlier paper \cite{DaMi} and references therein.  Among other topics of recent interest there are various types of limit theorems, see e.g.\ \cite{Bur1, Bur2, Erh}, large deviations probabilities, see e.g.\ \cite{Bur3, Kol}, and recurrence/transience of multivariate
AR models \cite{Zer}.

In this paper, we consider an autoregressive sequence $X_{t+1} =\alpha_{t+1} X_t+
\beta \xi_{t+1}$, $t=0,1,\ldots$ in $d$-dimensional Euclidean space with i.i.d.\ autoregressive coefficients $\{\alpha_t\}_{t\in\N}$, an independent i.i.d.\ noise sequence $\{\xi_t\}_{t\in\N}$,  constant noise scaling parameter
$\beta$ and initial value $X_0=0$, which is denoted as AR(1) process. Assuming that $\{\alpha_t\}_{t\in\N}$
are a.s.\ positive with $\E \log \alpha_t <0$ and that $\{\xi_t\}_{t\in\N}$ have
a finite first moment, the distributions of the elements of the autoregressive sequence converge
weakly to the unique limiting/stationary distribution. 

We modify the autoregressive sequence by introducing a certain ``restart'' mechanism that seems to be new in this setting. 
Next, we take  a series of autoregressive sequences with various coefficients and parameters, and consider limits for their stationary distributions and for the stationary distributions of their modifications, by assuming that the autoregressive coefficients tend to 1 and the noise scaling parameters tend to zero. When they converge to their limits with the proper speed, the stationary  distributions of the original sequences without restart mechanisms converge to a normal distribution. Introducing the restart mechanism makes the problem much richer and leads to a large class of limiting distributions.

The main advance of this paper is that we have managed to describe and analyse this class of limiting distributions, with obtaining
exact expressions for their densities and characteristic functions.
In particular, we provide conditions for the limiting distribution to be normal, which may be understood as conditions for the ``goodness'' of the restart mechanism. 
The obtained class of distributions has a number of interesting properties and seems to be new in the literature. 

In more detail, our ``restart'' mechanism works as follows.
We take a neighbourhood $A$ of the origin such that the sequence $\{X_t\}_{t\in\N_0}$ 
visits it with probability one. Then we run another autoregressive
sequence of the form $Y_{t+1}=\alpha_{t+1}Y_t{\bf 1}(Y_t\notin A)+\beta \xi_{t+1}$ with initial value $Y_0=0$, which restarts from the origin after each visit to $A$ (here ${\mathbf 1}$ is the indicator function). The modified sequence $\{Y_t\}_{t\in\N_0}$ has clearly
a regenerative structure. 
 
Further, we take a parametric family 
of multi-dimensional  autoregressive models.
We consider a limiting regime where the $\alpha$'s  and the $\beta$'s
depend on a parameter $m\in\N$ such that, for a fixed $a\in(0,\infty)$,
$1-\E \alpha_{m,t}\sim a/m\to0$, $1-\E \alpha_{m,t}^2\sim 2a/m\to0$ and $\beta_m \sim 1/\sqrt{m} \to 0$  as $m\to\infty$. 
We remark that  the stationary distributions of $\{X_t\}_{t\in\N_0}\equiv \{X_t^{(m)}\}_{t\in\N_0}$ converge to a normal distribution. 

For the sequence with restart mechanism we assume that the neighbourhood also depends on the parameter $m$ as 
$A^{(m)} = \gamma_m A$ where $A$ is a fixed neighbourhood of the origin and $\{\gamma_m\}_{m\in\N}$ are positive numbers.  We then show that,
depending on the behaviour of $\gamma_m$ as $m\to\infty$, 
 the stationary distributions
of the Markov chains $\{Y_t\}_{t\in\N_0}=\{Y_t^{(m)}\}_{t\in\N_0}$ have a variety of weak limits.
In particular, if $\{\gamma_m\}_{m\in\N}$ tends to zero very fast, we get the same limit as in the non-truncated 
case. In general, the limit is such that any of its one-dimensional projections on a line crossing the origin is 
a mixture of three distributions, an atomic distribution at the origin and two  distributions of signed absolute values of normal 
random variables. The whole collection of the limiting distributions is detailed in Theorem \ref{thm:LimitMixing}. 
Sufficient conditions for the limiting distribution to be continuous are provided in Theorem \ref{thm:Degenerate}.

The case that $\alpha_{m,t}$ is deterministic and that $1-\alpha_{m,t}\sim a/m$ appears in many different areas and is referred to as "heavy traffic" scaling. AR processes with these autoregressive coefficients are known as "almost non-stationary" AR processes and one is interested in estimating $a$ as $m\to\infty$ (see e.g.\ \cite{DaMi,Rah}). A process with autoregressive coefficient $1-a/m$ is obtained by applying the Euler-Maruyama method with step width $1/m$ to the Langevin stochastic differential equation
$$
dZ(s) = -a Z(s) ds + dB(s), \quad s\geq 0,
$$
with $Z(0)=0$, where $(B(s))_{s\geq0}$ is a standard Brownian motion. Heavy-traffic regimes for reflected AR(1) sequences were considered in \cite{BMR} (see also the references therein). For heavy-traffic regimes in the context of stochastic networks we refer to  e.g.\ \cite{Asm,Bor0}.

We believe that we have opened a window to an interesting and challenging research topic. 
A similar restart mechanism may be introduced and analysed
in a much broader setting, 
for a general stable time-homogeneous Markov chain $\{X_t\}_{t\in\N_0}$ in $\R^d$ (or in a general Polish space).
It is known that $\{X_t\}_{t\in\N_0}$ can be represented as a {\it 
stochastic recursion} $X_{t+1}=f(X_t,\xi_{t+1})$ with
i.i.d.\ driving sequence $\{\xi_t\}_{t\in\N}$ and a measurable function $f$,
thanks to the famous Skorokhod representation from the 50's 
(see e.g.\ \cite{Kif} or \cite{Bor}). We may introduce again 
the restart mechanism, with taking a neighbourhood  $A$ of the origin such that the Markov chain visits it with probability one. Then we run another Markov chain $Y_{t+1} = f(Y_t {\mathbf 1}\{Y_t \notin A\},\xi_{t+1})$ with initial value $Y_0=0$.  This Markov chain is regenerative again. The regenerative cycle length may have a sufficiently light-tailed distribution, and we may expect that if the set $A$ is {\it ``relatively small''} or/and is ``reasonably shaped'', the stationary distribution of the new Markov chain (which is the averaged distribution over the regenerative cycle) may be close to the stationary distribution of the original Markov chain.  However, if the set $A$ is {\it ``relatively big''}, then the sequence $\{Y_t\}_{t\in\N_0}$ may look similar to the i.i.d.\ sequence $\{f(0,\xi_t)\}_{t\in\N_0}$, so its stationary distribution may be far away from the stationary distribution of $\{X_t\}_{t\in\N_0}$. In this general setting, one
may address again the natural questions, {\it how to justify the meanings of ``relatively big/small''} and {\it 
what are other distributions that may appear in the intermediate
case}.

Our paper is organised as follows. In Section 2, we introduce formally the model and present Proposition \ref{prop11} (that summarises basic properties of autoregressive models) and  our main results, Theorems \ref{thm:LimitMixing} and \ref{thm:Degenerate}, complemented by a number of examples and comments. We also recall in Proposition \ref {thm:NoTruncation} a corresponding limiting result for sequences $\{X_t\}_{t\in\N_0}$. Then in Section 3 we give the proof of Theorem \ref{thm:Degenerate} and in Section 4 the proof of Theorem \ref{thm:LimitMixing}. 
To make the paper self-contained, we provide in Appendix \ref{AppendixA} arguments of the proof of Proposition \ref{prop11} and, in particular, two lemmas that include elements of originality and may be of their own interest. In Appendix \ref{sec:SecondProof} we prove some of our main results under stronger assumptions via the method of moment. We believe that this alternative approach gives some further insights.

\section{Basic model and main results}

We consider a series of autoregressive sequences $X_{t+1}= \alpha_{t+1} X_t + \beta \xi_{t+1}$ in $\R^d$ with particular choices of random autoregressive coefficients and of the noise scaling parameters,   and their truncated versions.

Let $a\in(0,\infty)$ be fixed. We assume that $\{\alpha_m\}_{m\in\N}$ are {\it positive} random 
variables with {\it finite first two moments} and $\{\beta_m\}_{m\in\N}$ are constants such that, as $m\to\infty$,
\begin{align}\label{alphabeta}
1-\mathbb{E} \alpha_m \sim \frac{a}{m} \quad  \mbox{and} \quad 1-\mathbb{E} \alpha_m^2 \sim \frac{2a}{m} \quad \mbox{as well as} \quad \beta_m \sim \frac{1}{\sqrt{m}}.
\end{align}
A comment on a possible extension to random $\beta$'s is given in Remark \ref{Rem5}.
We let $\{\xi_t\}_{t\in\N}$ be i.i.d.\ copies of a random vector $\xi$ in $\R^d$. For $m\in\N$, let $\{\alpha_{m,t}\}_{t\in\N}$ be i.i.d.\ copies of $\alpha_m$. Then let $X_{t+1}^{(m)} = \alpha_{m,t+1} X_t^{(m)}+\beta_m \xi_{t+1}$, $t=0,1,\ldots$ with $X_0^{(m)}=0$.

Let $\{\gamma_m\}_{m\in\N}$ be a sequence of positive real-valued numbers and let $A\subseteq\R^d$ be a measurable set such that
\begin{align*}
B^d(0,\underline{r})\subseteq A \subseteq B^d(0,\overline{r})
\end{align*}
with some $\underline{r},\overline{r}\in(0,\infty)$, where $B^d(x,r)$ stands for the $d$-dimensional open ball with radius $r$ and centre $x$. Let a sequence $\{Y_t^{(m)}\}_{t\in\N_0}$, $m\in\N$, be given by
\begin{align}\label{Ydef}
Y^{(m)}_{t+1}=\alpha_{m,t+1} Y_t^{(m)} \mathbf{1}\{Y_t^{(m)}\notin \gamma_m A\}+\beta_m \xi_{t+1}, \quad t\in\N_0, \quad \text{and} \quad Y_0^{(m)}=0. 
\end{align}
One can see that the elements of the sequences $\{X_t^{(m)}\}_{t\in\N_0}$ and $\{Y_t^{(m)}\}_{t\in\N_0}$ coincide until time
\begin{align*}
\tau^{(m)} = \min \{t\ge 1: \  X_t^{(m)}\in \gamma_m A\}\le\infty
\end{align*}
and then $\{Y_t^{(m)}\}_{t\in\N_0}$ restarts from the origin.

By $Y^{(m)}$ we denote a random vector having the stationary distribution of $\{Y_t^{(m)}\}_{t\in\N_0}$, i.e., 
\begin{equation}\label{eqn:DefinitionYm}
Y^{(m)}\overset{d}{=}\alpha_{m} Y^{(m)} \mathbf{1}\{Y^{(m)}\notin \gamma_m A\}+\beta_m \xi,
\end{equation}
where $\overset{d}{=}$ stands for equality in distribution and $\alpha_m$, $Y^{(m)}$ and $\xi$ are assumed to be independent. The aim of this paper is to analyse the asymptotic behaviour of the distributions of $Y^{(m)}$ as $m$ grows to infinity.

First, we formulate Proposition \ref{prop11} that summarises mostly known facts
on existence and uniqueness of the stationary distribution
and finiteness of corresponding moments. To make the paper self-contained, we provide comments related to the proof of Proposition \ref{prop11}
in the appendix.

We denote by $\|\cdot\|$ the Euclidean norm in ${\mathbb R}^d$.

\begin{prop}\label{prop11}
Let $m\in\N$ be such that $\mathbb{P}(\alpha_m\in[d/(d+1),1))>0$ and $\E\alpha_m<1$ and assume that $\E\xi=0$, $\E \|\xi\|^2<\infty$ and the matrix $\E\xi\xi^T$ has full rank.

(i) Then  
\begin{align*}
\tau^{(m)}<\infty \quad \mbox{a.s.}
\end{align*}
 and, moreover,
\begin{align*}
{\mathbb E} e^{c\tau^{(m)}}<\infty \quad
\mbox{for some}\quad c=c(m)>0.
\end{align*}

(ii) The unique stationary distribution of the Markov chain \eqref{Ydef} is given by
\begin{align*}
\pi^{(m)}(\cdot)=\frac{1}{\E \tau^{(m)}} \E \sum_{j=1}^{\tau^{(m)}} \mathbf{1}\{X_j^{(m)}\in\cdot\} \equiv 
\frac{1}{\E \tau^{(m)}} \E \sum_{j=1}^\infty \mathbf{1}\{X_j^{(m)}\in\cdot\} \mathbf{1}\{\tau^{(m)}\geq j\}
\end{align*}
and the distributions of $Y_t^{(m)}$ converge to it in the total variation norm as $t\to\infty$.

(iii) For any $c >1$, if ${\mathbb E} \alpha_m^c<1$ and ${\mathbb E} ||\xi||^{c}<\infty$,
then ${\mathbb E} ||Y^{(m)}||^{c}<\infty$, too.
\end{prop}

\begin{rema}\label{rema0} It follows from \eqref{alphabeta} that $\mathbb{E}\alpha_m = 1-\frac{a}{m}(1+o(1))$ and $m\Var(\alpha_m)\to0$ as $m\to\infty$. Consequently, for $m$ sufficiently large, the assumptions $\mathbb{P}(\alpha_m\in [d/(d+1),1))>0$ and $\E\alpha_m<1$ are satisfied and, thus, $Y^{(m)}$ is well defined. Since we study the asymptotic behaviour of 
$Y^{(m)}$ as $m\to\infty$,
this is sufficient for our purposes. Nevertheless, we tacitly assume throughout this paper that the random variables $(\alpha_m)_{m\in\N}$ satisfy $\mathbb{P}(\alpha_m\in[d/(d+1),1))>0$ and $\E\alpha_m<1$ for all $m\in\N$ so that $Y^{(m)}$ is well-defined for all $m\in\N$.
\end{rema}
   
We formulate now our main result. It says that the limiting distribution of $Y^{(m)}$ is a mixture of an atom at $0$ and an absolutely continuous distribution. For this random vector we provide the characteristic function, the density and the distributions of one-dimensional projections. 

By $\overset{d}{\longrightarrow}$ and $\overset{\mathbb{P}}{\longrightarrow}$ we denote convergence in distribution and convergence in probability, respectively.

\begin{theo}\label{thm:LimitMixing}
Assume that \eqref{alphabeta} holds and that $\E\xi=0$, $\E\|\xi\|^2<\infty$ and the matrix $\Sigma:=\E \xi \xi^T$ has full rank. Moreover, let
\begin{align}\label{EY}
\E Y^{(m)}\to \mu\in\R^d \quad \text{as} \quad m\to\infty
\end{align}
and
\begin{align*}
\E \tau^{(m)} \to \widehat{\tau} \in [1,\infty] \quad \text{as} \quad m\to\infty.
\end{align*}
If $\widehat{\tau}\in[1,\infty)$, assume additionally that there exists a random variable $\tau$ such that $\tau^{(m)}\overset{\mathbb{P}}{\longrightarrow} \tau$ as $m\to\infty$.
Then
\begin{align}\label{eqn:LimitMixing}
Y^{(m)} \overset{d}{\longrightarrow} Y:=B_1\cdot Z \quad \text{as} \quad m\to\infty,
\end{align}
where $B_1$ and $Z$ are independent. The random variable $B_1$ takes values 0 and 1 with probabilities
\begin{equation}\label{eqn:Definition_p}
p:={\mathbb P}(B_1=1)=1-{\mathbb P}(B_1=0)= \begin{cases} 1-\E\tau/\widehat{\tau}, & \quad \widehat{\tau}\in[1,\infty),\\ 1, & \quad \widehat{\tau}=\infty, \end{cases}
\end{equation}
and the $d$-dimensional random vector $Z$ has an absolutely continuous distribution that is characterised by the following properties:

(i) The characteristic function of $Z$ is
\begin{equation}\label{eqn:Formula_phi}
\varphi_Z(u)= 
\bigg( 1  
+ \mathbf{i} \frac{ \sqrt{2a} \langle u,\mu\rangle}{p\sqrt{u^T\Sigma u}} \int_0^{\frac{\sqrt{u^T\Sigma u}}{\sqrt{2a}}} \exp\bigg(\frac{t^2}{2}\bigg) \, dt \bigg) \exp\bigg(-\frac{u^T\Sigma u}{4a}\bigg), \quad u\in\R^d.
\end{equation}
 
(ii) The density of $Z$ is given by
\begin{align}\label{densZ}
f_Z(x)  = \frac{\sqrt{a}^d}{\sqrt{\det(\Sigma)}\sqrt{\pi}^d } \exp(-ax^T\Sigma^{-1}x) 
+ \widetilde{f}_Z(x), \quad x\in {\mathbb R}^d,
\end{align}
where, for odd dimensions $d=1,3,5,\ldots$,
\begin{align}\label{densodd}
\widetilde{f}_Z(x) =  \frac{(-2)^{\frac{d-1}{2}} a^{\frac{d+2}{2}}}{\sqrt{\det(\Sigma)}\kappa_{d-1}(d-1)!! p} \langle \Sigma^{-1} \mu, x\rangle h^{((d-1)/2)}(ax^T\Sigma^{-1}x), \quad x\in\R^d,
\end{align}
and for even $d=2,4,\ldots$,
\begin{align}\label{denseve}
\widetilde{f}_Z(x)
=  \frac{(-2)^{\frac{d}{2}} a^{\frac{d+3}{2}}}{\kappa_{d}d!! p} \langle \Sigma^{-1} \mu, x\rangle \int_{-\infty}^\infty h^{(d/2)}(a(x^T\Sigma^{-1}x+z^2)) \, dz, \quad x\in\R^d.
\end{align}
Here $h(s)=e^{-s}/\sqrt{s}$ for $s>0$, $h^{(k)}$ is its $k$th derivative, $\kappa_d$ the volume of the $d$-dimensional unit ball, and $(2k)!!=2k\cdot 2(k-1)\cdot \ldots \cdot 4\cdot 2$ 
the double factorials, for $k\in\N$ (note that $\kappa_0=1$ and $0!!=1$).

(iii) For any $v\in\R^d$,
\begin{equation}\label{eqn:vProjection}
\langle v, Z \rangle \overset{d}{=} \frac{\sqrt{v^T\Sigma v}}{\sqrt{2a}} B_{2,v} |N|,
\end{equation}
where $N$ and $B_{2,v}$ are two independent random variables, with $N$ having the standard normal distribution and $B_{2,v}$ having a two-point distribution, 
\begin{align}\label{defB}
\mathbb{P}(B_{2,v}=1)=\frac{1}{2}+\frac{\sqrt{\pi a}\langle v,\mu \rangle}{2p\sqrt{v^T\Sigma v}} \quad \text{and } \quad  \mathbb{P}(B_{2,v}=-1)=\frac{1}{2}-\frac{\sqrt{\pi a}\langle v,\mu \rangle}{2p\sqrt{v^T\Sigma v}}.
\end{align}
\end{theo}

\begin{rema}\label{Rem1}
It follows that, in conditions of Theorem
\ref{thm:LimitMixing}, we always have $
\frac{\sqrt{\pi a}\langle v,\mu \rangle}{p\sqrt{v^T\Sigma v}}\le 1$, so \eqref{defB} defines a probability distribution.  
\end{rema}

\begin{rema}\label{Rem2}
Note that $Z$ has a multivariate normal distribution if and only if $\mu = 0$. The latter condition holds if, say, the distribution of 
$\xi$ is symmetric and the set $A$ is symmetric too (e.g.\ a ball). From \eqref{eqn:DefinitionYm} and the independence of $\alpha_m$, $Y^{(m)}$ and $\xi$, we obtain
$$
\E Y^{(m)} = \E (\alpha_m Y^{(m)} \mathbf{1}\{Y^{(m)}\notin \gamma_m A\}+\beta_m \xi) =  \E\alpha_m \E Y^{(m)} - \E\alpha_m \E Y^{(m)} \mathbf{1}\{Y^{(m)}\in \gamma_m A\} 
$$
so that
\begin{equation}\label{eqn:ExpectationYm}
\E Y^{(m)} = -\frac{\E\alpha_m }{1- \E\alpha_m } \E Y^{(m)} \mathbf{1}\{Y^{(m)}\in \gamma_m A\}.
\end{equation}
This implies, together with $|\E Y^{(m)} \mathbf{1}\{Y^{(m)}\in \gamma_m A\}|\leq \overline{r} \gamma_m$ and \eqref{alphabeta}, that we have $\mu = \lim_{m\to\infty} \E Y^{(m)} = 0$ if
$m\gamma_m\to 0$ as $m\to\infty$.
\end{rema}

\begin{rema}\label{Rem6}
We get $\mu=0$ again if assume that $\gamma_m\to\infty$ and  the conditions of Theorem \ref{thm:LimitMixing} hold. 
Indeed, by \eqref{eqn:ExpectationYm} and Proposition \ref{prop11} (ii),
\begin{align}\label{eq:Rem4}
{\mathbb E} Y^{(m)}= - \frac{{\mathbb E} \alpha_m}{1-{\mathbb E} \alpha_m} \cdot \frac{1}{{\mathbb E} \tau^{(m)}}
\left( {\mathbb E} \beta_m\xi_1 \mathbf{1} \{\tau^{(m)} =1\} + 
{\mathbb E} Y^{(m)}_{\tau^{(m)}} \mathbf{1} \{\tau^{(m)}>1\}\right).
\end{align} 
Here $\frac{{\mathbb E} \alpha_m}{1-{\mathbb E} \alpha_m} = O(m)$. Since ${\mathbb E} \xi=0$, we get that
${\mathbb E} \xi_1 \mathbf{1}\{\tau^{(m)}=1\} = 
- {\mathbb E} \xi_1 \mathbf{1}\{\tau^{(m)}>1\}$, so 
the Euclidean norm of the first term in the parentheses in the right-hand side of \eqref{eq:Rem4} is 
\begin{align*}
\beta_m||{\mathbb E} \xi_1 \mathbf{1} \{ \beta_m\xi_1 \notin \gamma_m A\} ||
\le 
\beta_m{\mathbb E} ||\xi_1|| \mathbf{1} \{ ||\xi_1|| \ge \gamma_m \underline{r}/\beta_m\} 
\le \frac{\beta_m^2{\mathbb E}|| \xi_1||^2}{\gamma_m\underline{r}} = o(1/m).
\end{align*}
Since ${\mathbb E} \tau^{(m)}\ge 1$ and since the norm of the second term in the parentheses in the right-hand side of \eqref{eq:Rem4} does not exceed
\begin{align*}
\gamma_m\overline{r} {\mathbb P} (\tau^{(m)}>1) \le  
\gamma_m \overline{r} {\mathbb P} (||\xi_1||> \gamma_m\underline{r}/{\beta_m})
\le \frac{\beta_m^2\overline{r}{\mathbb E}||\xi_1||^2}{\gamma_m\underline{r}^2} = o(1/m),
\end{align*}
it follows that $\mu =\lim_{m\to\infty} {\mathbb E} Y^{(m)}=0$. We have used finiteness of the second 
moment of $\xi$ and the standard Markov inequality.\\
We may get $\mu=0$ under weaker assumptions on the
$\gamma$'s if allow the tail of $||\xi ||$ to have finite moments 
of a higher order. If, for example, 
${\mathbf E} ||\xi||^K <\infty$ for some $K>2$, then we may use the Markov inequality for the $K$'th moment to conclude that $\mu=0$ if $\gamma_m \cdot m^{\frac{K-2}{2(K-1)}}\to\infty$.
Further, if the exponential moment
 ${\mathbb E} \exp (c||\xi||)$ is finite for a sufficiently large  $c>0$, then the exponential Markov inequality gives that 
 $\mu=0$ if  $\gamma_m\cdot \frac{\sqrt{m}}{\log m} \to\infty$.
\end{rema}

\begin{rema}\label{Rem3}
It follows that, for the limiting vector $Z=(Z_1,\ldots,Z_d)$ in Theorem \ref{thm:LimitMixing}, the absolute value $|Z_i|$ of every of its coordinates has the same distribution as the absolute value of a normal random variable.
In addition, the marginal distribution of any projection of $Z$ on a direction orthogonal to $\mu$ is normal, and the distribution of a non-orthogonal projection is a mixture of distributions of signed absolute values of normal random variables. Therefore, 
it would look plausible for the limiting vector $Z$ to  
coincide in distribution with a random vector $(\psi_1N_1,\ldots,\psi_dN_d)$
where $(N_1,\ldots,N_d)$ is a multivariate normal vector and $(\psi_1,\ldots,\psi_d)$ an independent random vector whose coordinates take values $\pm 1$ only. However, as it follows from \eqref{eqn:Formula_phi}--\eqref{denseve},
this is not the case if $d>1$ and $\mu \ne 0$.
\end{rema}

\begin{rema}\label{Rem4}
The assumption that the matrix $\Sigma=\E \xi \xi^T$ has full rank is not restrictive. If $\Sigma$ is not regular, the components of $\xi$ and, thus, 
the components of $\{X_t^{(m)}\}_{t\in\N_0}$ are linearly dependent. In this case, it is sufficient to study a maximal subset of linearly independent components, for which the assumption on the covariance matrix takes place. 
\end{rema}

\begin{rema}\label{Rem5}
One may consider a more general case where the deterministic noise scaling parameter $\beta_m$ is replaced by i.i.d.\ random
variables $\{\beta_{m,t}\}_{t\in\N}$, for each $m$, and where $\E \beta_{m,1}\sim K/\sqrt{m}$ for a constant $K$,
as $m\to\infty$. This case may be reduced
to the case of constant $\beta$'s, by
introducing new random vectors $\xi_{m,t} = \sqrt{m} \beta_{m,t} \xi_t$ and new $\widehat{\beta}_m \sim 1/\sqrt{m}$. To find conditions on the $\xi_{m,t}$'s under which our results continue to hold,  
is a technical problem we do not address in this paper.  
\end{rema} 

\begin{rema}\label{rem:Subsequence}
Theorem \ref{thm:LimitMixing} can be applied to subsequences, i.e., if the assumptions of the theorem hold true for a subsequence $\{m_n\}_{n\in\N}$, then the statements stay valid along that subsequence.
Further, if 
all conditions of Theorem \ref{thm:LimitMixing} are satisfied except of \eqref{EY}, then Proposition \ref{prop:2ndMoment} 
implies convergence of  
the second moments of $\{Y^{(m)}\}_{m\in\N}$. 
In turn, this
yields tightness of the distributions of $\{Y^{(m)}\}_{m\in\N}$ and, for each convergent subsequence,  
convergence of the means. Therefore,  
one can apply Theorem \ref{thm:LimitMixing} to this subsequence. 
In particular, if 
$d=1$ and one knows additionally that $\lim_{m\to\infty} \mathbb{P}(Y^{(m)}\leq -y)=0$ for all $y\in(0,\infty)$, then 
each convergent subsequence must converge to the limiting distribution, say $Q$, with $\mu=\frac{p\sqrt{\E\xi^2}}{\sqrt{\pi a}}$. This implies convergence in distribution to $Q$ 
and convergence of the means too, 
$\E Y^{(m)}\to \frac{p\sqrt{\E\xi^2}}{\sqrt{\pi a}}$ as $m\to\infty$.    
\end{rema}

Sufficient conditions for the limiting distribution in Theorem \ref{thm:LimitMixing} to be continuous are given
in the following theorem. For $A_1,A_2\subset\R^d$ we define $A_1+A_2=\{x_1+x_2: x_1\in A_1, x_2\in A_2\}$.

\begin{theo}\label{thm:Degenerate}
Assume 
that $\alpha_m \overset{\mathbb{P}}{\longrightarrow} 1$ as $m\to\infty$ and that
\begin{align*}
\frac{\gamma_m}{\beta_m}\to\varrho\in[0,\infty) \quad \text{as} \quad m\to\infty.
\end{align*}
If there exists an $\varepsilon>0$ such that 
\begin{equation}\label{eqn:assumption_RW}
\E \inf\big\{t\geq 1: \sum_{i=1}^t \xi_i \in \varrho A+B^d(0,\varepsilon)\big\}=\infty, 
\end{equation}
then 
\begin{equation}\label{eqn:LimInf}
\liminf_{m\to\infty} \mathbb{P}(\tau^{(m)}>j)>0
\end{equation}
for $j\in\N$ and
\begin{align}\label{eqn:tauinfty}
\widehat{\tau}\equiv \lim_{m\to\infty} \E \tau^{(m)}=\infty.
\end{align}
\end{theo}
\begin{rema}
In particular, \eqref{eqn:assumption_RW} is satisfied if ${\mathbb E} \xi = 0$ and
if there exists a convex set $D\supseteq A$  such that 
\begin{align}\label{notin}
{\mathbb P} (\xi \notin \rho D + B^d(0,\varepsilon))>0,
\end{align} 
for some $\varepsilon>0$. The latter condition clearly holds if $\|\xi\|$ has unbounded support.
\end{rema}

\noindent {\bf Examples.} 
Consider a few simple examples. Let $d=1$ and $B_2=B_{2,1}$.\\
(1) 
Assume that  $\gamma_m=\beta_m=1/\sqrt{m}$, $\alpha_m=1-a/m$ and that
$\xi$ is uniformly distributed in the interval $(-1,1)$. \\
(1.1)  Assume first that $A=(-1/2,1/2)$.  
 Then ${\mathbb P} (\tau^{(m)}\ge 2) = 1/2$ and
 \eqref{notin} holds.  So, by Theorem \ref{thm:Degenerate}, ${\mathbb E} \tau^{(m)} \to\infty$ as $m\to\infty$. By symmetry, we have $\mu=0$, and the conclusion of Theorem \ref{thm:LimitMixing} holds with ${\mathbb P} (B_1=1)=1$ and ${\mathbb P}(B_2=1)={\mathbb P}(B_2=-1)=1/2$.\\
 (1.2) If instead $A=(-1,1/2)$, then, by Theorem \ref{thm:Degenerate}, ${\mathbb E} \tau^{(m)} \to\infty$ as $m\to\infty$. Since ${\mathbb P} (Y_t^{(m)}\ge -1/\sqrt{m})=1$ for all $t,m$, as discussed in Remark \ref{rem:Subsequence}, Theorem \ref{thm:LimitMixing} yields ${\mathbb P} (B_1=1)=1$ and  ${\mathbb P}(B_2=1)=1$.\\ 
 (2) Assume that  $\gamma_m=\beta_m=1/\sqrt{m}$, that $\xi$ is uniformly distributed in the interval $(-1,1)$ and that $\alpha_{m,t}=\widetilde{\alpha}_t -\frac{1}{2m}$ where $\widetilde{\alpha},\widetilde{\alpha}_1, \widetilde{\alpha}_2,\ldots$ are i.i.d random variables with ${\mathbb P} (\widetilde{\alpha}> 1/2)=1$, ${\mathbb E} \widetilde{\alpha}  =1$ and ${\mathbb E} \log \widetilde{\alpha} <0$
(for example, we may assume that ${\mathbb P}(\widetilde{\alpha} =5/4)={\mathbb P}(\widetilde{\alpha}  =3/4)=1/2)$. 
\\
Then, for any $m$ and $t$, we have $|Y_t^{(m)}|\le \frac{Z_t}{\sqrt{m}}$ where
\begin{align*}
Z_0=0 \ \mbox{and} \ Z_{t+1}=\widetilde{\alpha}_{t+1} Z_t + |\xi_{t+1}|, \ \mbox{for} \ t\ge 0.
\end{align*}
It is well-known that, for any $c>0$,
the hitting time $\tau_c^Z = \min \{t>0: Z_t\le c\}$ has a finite first moment. Therefore,
$\widehat{\tau}<\infty$ for any choice of neighbourhood $A$ of zero. Further, by the Lebesgue theorem,
$\{\tau^{(m)}\}_{m\in\N}$ are uniformly integrable, so $Y^{(m)}\to0$ in probability and in $\mathcal{L}^1$ as $m\to\infty$.

The following proposition provides an example showing that, in the case $d=1$, the
sequence ${\mathbb E} Y^{(m)}$ may converge to any
$\mu$ from a non-degenerate interval.

\begin{prop}\label{Prop2.4}
Let $d=1$.
Let $\alpha_m=1-a/m$ for $m\in\N$ and let $\xi$ be such that $\E\xi=0$, $\E \xi^2=1$ and the distribution of $\xi$ has a density. Assume that
$$
A=[-2,1] \text{ and $\xi$ has support $[-1,1]$}
$$
or that 
$$
A=[-1,1] \text{ and $\xi$ has support $[-1/2,1]$}.
$$
Then, for any $\mu\in [0,1/\sqrt{\pi a}]$, there exists a sequence $\{\gamma_m\}_{m\in\N}$ such that
$$
\E Y^{(m)}\to\mu \quad \text{and} \quad \E\tau^{(m)}\to\infty \quad \text{as} \quad m\to\infty.
$$ 
\end{prop}

\begin{rema}
Let $N$ be a standard normal random variable, $F_1$ the distribution of $|N|/\sqrt{2a}$, and 
$F_2$ the distribution of $-|N|/\sqrt{2a}$. Under the assumptions of Proposition \ref{Prop2.4}, for any $c\in[1/2,1]$ one can choose $\{\gamma_m\}_{m\in\N}$ such that the limiting distribution in Theorem \ref{thm:LimitMixing} is the mixture $cF_1+(1-c)F_2$. By replacing $A$ and $\xi$ by $-A$ and $-\xi$ one gets the corresponding result for $c\in[0,1/2]$.
\end{rema}

To complete the section, we make a few comments about the autoregressive sequences $\{X^{(m)}_t\}_{t\in\N_0}$, $m\in\N$. By $X^{(m)}$ we denote a random vector with the stationary distribution of $\{X_t^{(m)}\}_{t\in\N_0}$, i.e.,
$$
X^{(m)}\overset{d}{=}\alpha_m X^{(m)} + \beta_m\xi,
$$
where $\alpha_m$, $\xi$ and $X^{(m)}$ are assumed to be independent. If ${\mathbb E} \log \alpha_m <0$, the stationary distribution uniquely exists and $X_t^{(m)}\overset{d}{\longrightarrow}X^{(m)}$ as $t\to\infty$. This follows from Theorem 1.6 and Theorem 1.5 in \cite{Vervaat1979}. For the sequence $\{X^{(m)}\}_{m\in\N}$ we have the following limit theorem. 

\begin{prop}\label{thm:NoTruncation}
Assume that \eqref{alphabeta} is satisfied as well as $\E\xi=0$ and $\E\|\xi\|^2<\infty$. Let $N_{\widehat{\Sigma}}$ be a centred normal random vector with covariance matrix $\widehat{\Sigma}:=\E\xi\xi^T/(2a)$. Then, $X^{(m)}\overset{d}{\longrightarrow} N_{\widehat{\Sigma}}$ as $m\to\infty$.
\end{prop}

This result should be known too. However, we could not find a proper reference and, to make the paper self-contained, we decided to comment on its proof -- see the end of Section 4.

\section{Proofs of Theorem \ref{thm:Degenerate} and Proposition \ref{Prop2.4}}

\begin{proof}[Proof of Theorem \ref{thm:Degenerate}.] 
We define a random walk 
$$
S_t = \sum_{i=1}^t \xi_i, \quad t\in\N, \quad \text{and} \quad S_0=0
$$
as well as random walks $S_t^{(m)}=\beta_m S_t$ for $m\in\N$ and $t\in\N_0$.
Let $\theta = \min \{t\ge 1: S_t\in \varrho A+B^d(0,\varepsilon) \}$ be the first hitting time of the set $\varrho A+B^d(0,\varepsilon)$ by the random walk $\{S_t\}_{t\in\N_0}$. By assumption \eqref{eqn:assumption_RW}, we have $\E \theta = \infty$. 

In what follows we will link $\{X_t^{(m)}\}_{t\in\N_0}$ with the  random walk $\{S_t^{(m)}\}_{t\in\N_0}$.
For $j\in\N$ and $m\in\N$ we have that 
\begin{align*}
& \mathbb{P}(\tau^{(m)}>j) \\
& \geq \mathbb{P}(X_\ell^{(m)} \notin \gamma_m A, \ell\in\{1,\hdots,j\}) \\
& \geq \mathbb{P}(S_\ell^{(m)}\notin \gamma_m A+B^d(0,\beta_m\varepsilon/2), \ell\in\{1,\hdots,j\})- \mathbb{P}(\max_{\ell\in\{1,\hdots,j\}} ||S_\ell^{(m)}-X_\ell^{(m)}||\geq \beta_m\varepsilon/2)\\
& \geq \mathbb{P}(S_\ell\notin \frac{\gamma_m}{\beta_m} A + B^d(0,\varepsilon/2), \ell\in\{1,\hdots,j\})- \mathbb{P}(\max_{\ell\in\{1,\hdots,j\}} ||S_\ell^{(m)}-X_\ell^{(m)}||\geq \beta_m\varepsilon/2).
\end{align*}
For $i,u\in\N_0$ with $i\le u$, let $\Lambda_{m,i,u} = \prod_{k=i+1}^u \alpha_{m,k}$ if $i<u$ and $\Lambda_{m,i,u}=1$ if $i=u$.
From $\alpha_m\overset{\mathbb{P}}{\longrightarrow} 1$ as $m\to\infty$ it follows that $\Lambda_{m,i,u}\overset{\mathbb{P}}{\longrightarrow} 1$ as $m\to\infty$.
Clearly, for $\ell\in\N$, we have
\begin{equation}\label{eqn:Xlm}
X_\ell^{(m)}=\sum_{i=1}^\ell \Lambda_{m,i,\ell} \beta_m \xi_i
\end{equation}
and
\begin{equation}\label{eqn:Slm-Xlm}
S_\ell^{(m)}-X_\ell^{(m)} = \sum_{i=1}^\ell (1-\Lambda_{m,i,\ell}) \beta_m \xi_i.
\end{equation}
Therefore, we obtain
$$
\mathbb{P}(||S_\ell^{(m)}-X_\ell^{(m)}||\geq \beta_m\varepsilon/2) \leq
\sum_{i=1}^\ell \mathbb{P} \bigg(|1-\Lambda_{m,i,l}| \, \|\xi_i\|>\frac{\varepsilon}{2\ell}\bigg),
$$
where each term in the right-hand side tends to zero as $m\to\infty$, 
so that
$$
\lim_{m\to\infty} \mathbb{P}(\max_{\ell\in\{1,\hdots,j\}} ||S_\ell^{(m)}-X_\ell^{(m)}||\geq \beta_m\varepsilon/2) =0.
$$
For $m$ sufficiently large we have $\frac{\gamma_m}{\beta_m} A + B^d(0,\varepsilon/2)\subseteq \varrho A+ B^d(0,\varepsilon)$ and, thus,
$$
\mathbb{P}(S_\ell\notin \frac{\gamma_m}{\beta_m} A + B^d(0,\varepsilon/2), \ell\in\{1,\hdots,j\}) \geq \mathbb{P}(S_\ell\notin \varrho A + B^d(0,\varepsilon), \ell\in\{1,\hdots,j\}).
$$
Thus, we have shown that
$$
\liminf_{m\to\infty}\mathbb{P}(\tau^{(m)}>j) \geq \mathbb{P}(S_\ell\notin \varrho A+B^d(0,\varepsilon), \ell\in\{1,\hdots,j\}) = \mathbb{P}(\theta>j),
$$
which is \eqref{eqn:LimInf}. Consequently, the Fatou's lemma yields
$$
\liminf_{m\to\infty} \E \tau^{(m)} = \liminf_{m\to\infty} \sum_{j=0}^\infty \mathbb{P}(\tau^{(m)}>j) \geq \sum_{j=0}^\infty \mathbb{P}(\theta>j) = \E \theta = \infty.
$$
This completes the proof of \eqref{eqn:tauinfty}.
\end{proof}

\begin{proof}[Proof of Proposition \ref{Prop2.4}]
We prove the statement simultaneously for both choices of the set $A$ and the support of $\xi$. Since the distribution of $\xi$ is absolutely continuous, the joint distribution of $(X^{(m)}_1,\hdots,X^{(m)}_{t})$ is absolutely continuous, for any $t\in\N$. This yields that the maps
$$
\gamma_m \mapsto \E \tau^{(m)} = \E\sum_{t=1}^{\infty} \mathbf{1}\{ X_1^{(m)},\hdots,X^{(m)}_{t-1}\notin \gamma_m A \}
$$
and
$$
\gamma_m \mapsto \E \sum_{t=1}^{\infty} X_t^{(m)} \mathbf{1}\{ X_1^{(m)},\hdots,X^{(m)}_{t-1}\notin \gamma_m A \}
$$
are continuous.
Indeed, by the dominated convergence theorem, one can interchange the limits in $\gamma_m$ with the expectations and the sums. 
Integrable upper bounds are obtained by taking a smaller $\gamma_m$. 
Combining the continuity of the two functions above, we can see that $\gamma_m\mapsto \E Y^{(m)}$ is continuous.

Due to \eqref{eqn:ExpectationYm}, $\E Y^{(m)}$ becomes very close to $0$ if we choose $\gamma_m\leq 3\beta_m/4$ sufficiently small.

Let $\gamma_m=c \beta_m$, $m\in\N$, for a fixed $c\in (1/2,1)$. Then we have $\mathbb{P}(Y^{(m)}\leq -\gamma_m)=0$. From Theorem \ref{thm:Degenerate} it follows that $\widehat{\tau}=\infty$. As described in Remark \ref{rem:Subsequence}, we have $\E Y^{(m)}\to 1/\sqrt{\pi a}$ as $m\to\infty$.

The above arguments show that, for $m$ sufficiently large, one can choose $\gamma_m\leq 3\beta_m/4$ such that $\E Y^{(m)}$ is close to zero or such that $\E Y^{(m)}$ is close to $1/\sqrt{\pi a}$. By the continuity of $\gamma_m\mapsto \E Y^{(m)}$ and the intermediate value theorem, $\gamma_m\leq 3\beta_m/4$ can be chosen in such a way that any value in between is attained. Due to $\gamma_m\leq 3\beta_m/4$, it follows from the fact that, by convexity of $A$, $\E\tau^{(m)}$ is decreasing in $\gamma_m/\beta_m$ and Theorem \ref{thm:Degenerate} that $\widehat{\tau}=\infty$ for this choice of $\{\gamma_m\}_{m\in\N}$.
\end{proof}

\section{{Proofs of Theorem \ref{thm:LimitMixing} and Proposition \ref{thm:NoTruncation}} }

This section includes four subsections, three of them (Subsections 4.1--4.3) are devoted
to the proof of Theorem \ref{thm:LimitMixing}, and the last
subsection provides a short proof of Proposition \ref{thm:NoTruncation}. 
 Recall that we assume tacitly that $\mathbb{P}(\alpha_m\in[d/(d+1),1))>0$ and $\E\alpha_m<1$ for all $m\in\N$, see Remark \ref{rema0}.

\subsection{Preliminaries}

\begin{lemm}\label{lem:BoundMomentInterval}
Assume that $\E \|\xi\|^k<\infty$ for some $k\in\N$. Then
$$
\E  |\langle u,Y^{(m)} \rangle|^k \mathbf{1}\{Y^{(m)}\in \gamma_m A\} \leq \|u\|^k \bigg(1 + \frac{\overline{r}^k}{\underline{r}^k} \bigg) \beta_m^k \frac{\E \|\xi\|^k}{\E \tau^{(m)}}
$$
for all $u\in\R^d$. 
\end{lemm}

\begin{proof}
It follows from Proposition \ref{prop11} (ii) and the elementary observation
\begin{equation}\label{eqn:BoundInnerProductA}
\max_{x\in A} |\langle u, x\rangle|\leq \max_{x\in B^d(0,\overline{r})} |\langle u, x\rangle|\leq \|u\| \overline{r}
\end{equation}
that
\begin{align*}
 \E |\langle u, Y^{(m)} \rangle|^k \mathbf{1}\{Y^{(m)}\in \gamma_m A\} 
& = \frac{1}{\E \tau^{(m)}} \E \sum_{j=1}^\infty \mathbf{1}\{\tau^{(m)}\geq j\} |\langle u,X_j^{(m)}\rangle|^k \mathbf{1}\{X_j^{(m)}\in \gamma_m A\} \allowdisplaybreaks\\
& = \frac{1}{\E \tau^{(m)}} \E \sum_{j=1}^\infty \mathbf{1}\{\tau^{(m)}= j\} |\langle u,X_j^{(m)}\rangle|^k \\
& \leq \frac{\E \mathbf{1}\{\beta_m \xi_1\in \gamma_m A\} \beta_m^k|\langle u,\xi_1\rangle|^k  + \mathbb{P}(\tau^{(m)}\geq 2) \gamma_m^k \|u\|^k \overline{r}^k }{\E \tau^{(m)}}.
\end{align*}
By the Markov inequality, we have that
$$
\mathbb{P}(\tau^{(m)}\geq 2) = \mathbb{P}(\beta_m \xi \notin \gamma_m A) \leq \mathbb{P}(\|\xi\|\geq \gamma_m\underline{r}/\beta_m) \leq \frac{\beta_m^k \E \|\xi\|^k}{\gamma_m^k \underline{r}^k}.
$$
Together with 
$$
\E \mathbf{1}\{\beta_m \xi_1\in \gamma_m A\} \beta_m^k|\langle u, \xi_1\rangle|^k \leq \beta_m^k \E |\langle u,\xi_1\rangle |^k \leq  
\beta_m^k \|u\|^k \E \|\xi\|^k,
$$
this concludes the proof.
\end{proof}

\begin{lemm}\label{lemm:2ndMoment}
Assume that $\E\xi=0$ and $\E\|\xi\|^2<\infty$ and that assumption \eqref{alphabeta} holds.

(i) If $\E\tau^{(m)}\to\infty$ as $m\to\infty$, then
$$
\lim_{m\to\infty} m \E \langle u, Y^{(m)}\rangle^2 \mathbf{1}\{Y^{(m)}\in\gamma_m A_m\}=0
$$
for all $u\in\R^d$.

(ii) If $\E\tau^{(m)}\to\widehat{\tau}\in[1,\infty)$ as $m\to\infty$ and if there exists a random variable $\tau$ on the same probability space such that $\tau^{(m)}\overset{\mathbb{P}}{\longrightarrow}\tau$ as $m\to\infty$, then
$$
\lim_{m\to\infty} m \E \langle u, Y^{(m)}\rangle^2 \mathbf{1}\{Y^{(m)}\in\gamma_m A_m\}= \frac{\E\tau}{\widehat{\tau}} \E \langle u, \xi\rangle^2
$$
for all $u\in\R^d$.
\end{lemm}

\begin{proof}
We start with part (i) of the lemma. From Lemma \ref{lem:BoundMomentInterval}, we obtain that
$$
m \E \langle u,Y^{(m)}\rangle^2 \mathbf{1}\{Y^{(m)}\in\gamma_m A\} \leq m\beta_m^2 \|u\|^2 \bigg(1+\frac{\overline{r}^2}{\underline{r}^2} \bigg) \frac{\E\|\xi\|^2}{\E \tau^{(m)}}. 
$$
Since, by assumption, $\E \tau^{(m)}\to\infty$ as $m\to\infty$, the right-hand side vanishes as $m\to\infty$. This proves part (i).

Next we consider part (ii). First we assume that $\mathbb{P}(\tau>1)=0$. It follows from Proposition \ref{prop11} (ii) and \eqref{eqn:BoundInnerProductA} that
\begin{align*}
& \bigg|m \E \langle u,Y^{(m)}\rangle^2 \mathbf{1}\{Y^{(m)}\in\gamma_m A\} - \frac{\E\langle u,\xi\rangle^2}{\E\tau^{(m)}} \bigg| \\
& = \frac{1}{\E \tau^{(m)}} \bigg| m\E\sum_{j=1}^\infty \mathbf{1}\{\tau^{(m)}=j\} \langle u,X_j^{(m)}\rangle^2 - \E\langle u,\xi\rangle^2\bigg| \\
& = \frac{1}{\E \tau^{(m)}} \bigg| m \E \langle u,X_1^{(m)}\rangle^2 - \E \langle u,\xi\rangle^2 - m\E \mathbf{1}\{\tau^{(m)}>1\} \langle u,X_1^{(m)}\rangle^2 \\
& \quad \quad \quad \quad \quad + m\E\sum_{j=2}^\infty \mathbf{1}\{\tau^{(m)}=j\} \langle u,X_j^{(m)}\rangle^2 \bigg| \\
& \leq \frac{\E\langle u,\xi\rangle^2}{\E\tau^{(m)}} |m\beta_m^2-1|  + \frac{m\beta_m^2}{\E \tau^{(m)}}  \E \langle u,\xi_1\rangle^2\mathbf{1}\{\tau^{(m)}>1\} + \frac{m\|u\|^2\overline{r}^2\gamma_m^2}{\E \tau^{(m)}} \mathbb{P}(\tau^{(m)}>1).
\end{align*}
We have
\begin{align*}
m\gamma_m^2 \mathbb{P}(\tau^{(m)}>1) 
&= m\gamma_m^2 \mathbb{P}(\beta_m\xi \notin \gamma_m A) \\
& \leq m\gamma_m^2 \E \frac{\beta_m^2\|\xi\|^2}{\gamma_m^2\underline{r}^2} \mathbf{1}\{\beta_m\xi\notin\gamma_m A\} = \frac{m\beta_m^2}{\underline{r}^2} \E \|\xi_1\|^2 \mathbf{1}\{\tau^{(m)}>1\}
\end{align*}
whence
\begin{align*}
& \big|m \E \langle u,Y^{(m)}\rangle^2 \mathbf{1}\{Y^{(m)}\in \gamma_m A\} - \frac{\E\langle u,\xi\rangle^2}{\E\tau^{(m)}} \big|\\
& \leq \frac{\E\langle u,\xi\rangle^2}{\E\tau^{(m)}} |m\beta_m^2-1| + \bigg(1+ \frac{\overline{r}^2}{\underline{r}^2}\bigg) \|u\|^2 \frac{m\beta_m^2}{\E \tau^{(m)}}  \E \|\xi_1\|^2\mathbf{1}\{ \tau^{(m)}>1\}.
\end{align*}
The assumption $\mathbb{P}(\tau>1)=0$ implies $\mathbb{P}(\tau^{(m)}>1)\to 0$ as $m\to\infty$ and, thus,
$$
\lim_{m\to\infty} \E \|\xi_1\|^2\mathbf{1}\{ \tau^{(m)}>1\} =0.
$$
Together with $m\beta_m^2\to1$ and $\E\tau^{(m)}\to\widehat{\tau}$ as $m\to\infty$, we obtain
$$
\lim_{m\to\infty} m \E \langle u, Y^{(m)}\rangle^2 \mathbf{1}\{Y^{(m)}\in\gamma_m A\} = \frac{\E\langle u,\xi\rangle^2}{ \widehat{\tau}},
$$
which completes the proof of (ii) for $\mathbb{P}(\tau>1)=0$.

Next we assume $\mathbb{P}(\tau>1)>0$. Note that
$$
\mathbb{P}(\tau^{(m)}>1) = \mathbb{P}(\beta_m \xi\notin\gamma_m A) \leq \mathbb{P}(\beta_m\|\xi\|\geq \gamma_m \underline{r}).
$$
For a subsequence of $\gamma_m/\beta_m$ that converges to infinity, the right-hand side vanishes, which contradicts $\mathbb{P}(\tau>1)>0$. Thus, we obtain $\varrho:=\limsup_{m\to\infty} \gamma_m/\beta_m<\infty$. It follows from Proposition \ref{prop11} (ii) that
$$
m \E \langle u,Y^{(m)}\rangle^2 \mathbf{1}\{Y^{(m)}\in\gamma_m A\} 
= \frac{m}{\E \tau^{(m)}} \E\sum_{j=1}^\infty \mathbf{1}\{\tau^{(m)}=j\} \langle u,X_j^{(m)}\rangle^2.
$$
For $k\in\N$ we can decompose the right-hand side into the sum of
$$
T_{1,k,m}:= \frac{m}{\E \tau^{(m)}} \E\sum_{j=1}^k \mathbf{1}\{\tau^{(m)}=j\} \langle u,X_j^{(m)}\rangle^2
$$
and
$$
T_{2,k,m}:= \frac{m}{\E \tau^{(m)}} \E\sum_{j=k+1}^\infty \mathbf{1}\{\tau^{(m)}=j\} \langle u,X_j^{(m)}\rangle^2.
$$
We obtain that
\begin{align*}
T_{2,k,m} \leq \frac{m}{\E \tau^{(m)}} \E \mathbf{1}\{\tau^{(m)}\geq k+1\} \|u\|^2 \overline{r}^2 \gamma_m^2 
 \leq \|u\|^2 \overline{r}^2 \frac{m\gamma_m^2 \E \tau^{(m)}}{\E \tau^{(m)} (k+1)}, 
\end{align*}
where we used the Markov inequality in the second step. This implies that
\begin{equation}\label{eqn:LimitT2km}
\limsup_{m\to\infty} T_{2,k,m} \leq \frac{ \overline{r}^2 \varrho^2 \|u\|^2}{k+1}.
\end{equation}
As in the proof of Theorem \ref{thm:Degenerate} let $S_j^{(m)}=\beta_m \sum_{i=1}^j \xi_i$ and $S_j=\sum_{i=1}^j \xi_i$ for $j,m\in\N$. 
We have that
\begin{align*}
\widetilde{T}_{1,k,m} & := T_{1,k,m} - \frac{m}{\E \tau^{(m)}}  \E \sum_{j=1}^k \mathbf{1}\{\tau^{(m)}=j\} \langle u,S_j^{(m)}\rangle^2 \\
& = \frac{m}{\E \tau^{(m)}}  \E \sum_{j=1}^k \mathbf{1}\{\tau^{(m)}=j\} \big(\langle u,X_j^{(m)}\rangle^2 - \langle u,S_j^{(m)}\rangle^2\big).
\end{align*}
For any $j\in\N$, we obtain that
\begin{align*}
\big| \E \langle u,X_j^{(m)}\rangle^2 - \langle u,S_j^{(m)}\rangle^2\big| & \leq \sqrt{\E \langle u,S_j^{(m)} + X_j^{(m)} \rangle^2} \sqrt{\E \langle u,S_j^{(m)} - X_j^{(m)} \rangle^2} \\
& \leq \|u\|^2 \sqrt{\E \|S_j^{(m)} + X_j^{(m)} \|^2} \sqrt{\E \|S_j^{(m)} - X_j^{(m)} \|^2}.
\end{align*}
Combining \eqref{eqn:Xlm}  and \eqref{eqn:Slm-Xlm} with straightforward computations leads to
$$
\E \|S_j^{(m)} + X_j^{(m)} \|^2 = \beta_m^2 \sum_{i=1}^j (1 + 2 (\E \alpha_m)^{j-i} + (\E \alpha_m^2)^{j-i}) \E \|\xi_i\|^2
$$
and  
$$
\E \|S_j^{(m)} - X_j^{(m)} \|^2 = \beta_m^2 \sum_{i=1}^j (1 - 2 (\E \alpha_m)^{j-i} + (\E \alpha_m^2)^{j-i}) \E \|\xi_i\|^2.
$$
Since, by \eqref{alphabeta}, $m\beta_m^2\to1$, $(1 - 2 (\E \alpha_m)^{j-i} + (\E \alpha_m^2)^{j-i})\to 0$ and $(1 + 2 (\E \alpha_m)^{j-i} + (\E \alpha_m^2)^{j-i})\to 4$ as $m\to\infty$, we obtain
$$
\lim_{m\to\infty} |\widetilde{T}_{1,k,m}| =0.
$$
It follows from \eqref{alphabeta} and $\tau^{(m)}\overset{\mathbb{P}}{\longrightarrow}\tau$ as $m\to\infty$ that
$$
\lim_{m\to\infty} T_{1,k,m} = \lim_{m\to\infty} \frac{m}{\E \tau^{(m)}}  \E \sum_{j=1}^k \mathbf{1}\{\tau^{(m)}=j\} \langle u,S_j^{(m)}\rangle^2 = \frac{1}{\widehat{\tau}} \E \sum_{j=1}^k \mathbf{1}\{\tau=j\} \langle u,S_j\rangle^2.
$$
Letting $k\to\infty$ and using \eqref{eqn:LimitT2km} and the monotone convergence theorem, we deduce that
$$
\lim_{m\to\infty} m \E \langle u,Y^{(m)}\rangle^2 \mathbf{1}\{Y^{(m)}\in\gamma_m A\} = \frac{1}{\widehat{\tau}} \E \sum_{j=1}^\infty \mathbf{1}\{\tau=j\} \langle u,S_j\rangle^2 = \frac{1}{\widehat{\tau}} \E \langle u,S_{\tau}\rangle^2.
$$
For $j\in\N$ we have that $\mathbf{1}\{\tau^{(m)}\geq j\}$ is independent of $\xi_j$ for all $m\in\N$, whence $\mathbf{1}\{\tau\geq j\}$ is independent of $\xi_j$. Thus, by the Wald identity for the second moments (see \cite{Neveu}, pp.\ 72-74), we have $\E S^2_{\tau}=\E \tau \E \langle u,\xi\rangle^2$, which completes the proof.
\end{proof}

\begin{prop}\label{prop:2ndMoment}
Assume that $\E\xi=0$ and $\E\|\xi\|^2<\infty$ and that assumption \eqref{alphabeta} holds.

(1) If $\E\tau^{(m)}\to\infty$ as $m\to\infty$, then
$$
\lim_{m\to\infty} \E \langle u, Y^{(m)}\rangle^2=\frac{\E \langle u, \xi\rangle^2}{2a}
$$
for all $u\in\R^d$.

(2) If $\E\tau^{(m)}\to\widehat{\tau}\in[1,\infty)$ as $m\to\infty$ and if there exists a random variable $\tau$ on the same probability space such that $\tau^{(m)}\overset{\mathbb{P}}{\longrightarrow}\tau$ as $m\to\infty$, then
$$
\lim_{m\to\infty} \E \langle u, Y^{(m)}\rangle^2=\bigg( 1 - \frac{\E\tau}{\widehat{\tau}}\bigg)\frac{\E \langle u, \xi\rangle^2}{2a}
$$
for all $u\in\R^d$.
\end{prop}

\begin{proof}
It follows from \eqref{eqn:DefinitionYm} that
\begin{align*}
\E \langle u, Y^{(m)}\rangle^2 & = \E (\alpha_m \langle u, Y^{(m)} \rangle \mathbf{1}\{Y^{(m)}\notin \gamma_m A\}+\beta_m \langle u, \xi \rangle)^2 \\
& = \E\alpha_m^2 \E \langle u,Y^{(m)}\rangle^2 - \E\alpha_m^2 \E \langle u,Y^{(m)}\rangle^2 \mathbf{1}\{Y^{(m)}\in\gamma_m A\} + \beta_m^2\E\langle u,\xi\rangle^2,
\end{align*}
whence
\begin{equation}\label{eqn:Formula2ndMoment}
\E \langle u, Y^{(m)}\rangle^2 = \frac{1}{1-\E \alpha_m^2} \big( -\E \alpha_m^2 \E \langle u, Y^{(m)}\rangle^2 \mathbf{1}\{Y^{(m)}\in \gamma_m A\} + \beta_m^2 \E\langle u,\xi\rangle^2 \big).
\end{equation}
For the second term on the right-hand side \eqref{alphabeta} leads to
$$
\lim_{m\to\infty} \frac{\beta_m^2}{1-\E\alpha_m^2} \E\langle u,\xi\rangle^2=\frac{\E \langle u,\xi\rangle^2}{2a}.
$$
So it remains to study the asymptotic behaviour of the first term on the right-hand side of \eqref{eqn:Formula2ndMoment}. Since, by \eqref{alphabeta},
$$
\lim_{m\to\infty} \frac{\E\alpha_m^2}{m(1 - \E\alpha_m^2)}=\frac{1}{2a},
$$
we have
$$
 \lim_{m\to\infty} \frac{-\E \alpha_m^2 }{1-\E \alpha_m^2}  \E \langle u, Y^{(m)}\rangle^2 \mathbf{1}\{Y^{(m)}\in \gamma_m A\}  = -\frac{1}{2a}\lim_{m\to\infty} m \E \langle u,Y^{(m)}\rangle^2 \mathbf{1}\{Y^{(m)}\in\gamma_m A\}.
$$
Thus, the assertions follow from Lemma \ref{lemm:2ndMoment}.
\end{proof}

\begin{lemm}\label{lem:LimInd}
Let $u\in\mathbb{R}^d$. Assume that $\E\xi=0$ and $\E\|\xi\|^2<\infty$, that assumption \eqref{alphabeta} holds and that
$$
\E Y^{(m)}\to\mu\in\R^d \quad \text{and} \quad m \E \langle u, Y^{(m)} \rangle^2 \mathbf{1}\{Y^{(m)}\in\gamma_m A\}  \to w\in[0,\infty) \quad \text{as} \quad m\to\infty.
$$
Then,
$$
\lim_{m\to\infty} m \E \mathbf{1}\{Y^{(m)}\in\gamma_mA\} (1-e^{\mathbf{i} \alpha_m \langle u,Y^{(m)}\rangle})=\mathbf{i} a\langle u,\mu\rangle + \frac{w}{2}.
$$
\end{lemm}

\begin{proof}
By \eqref{eqn:DefinitionYm} we have
\begin{align*}
\E \langle u,Y^{(m)}\rangle & = \E \alpha_m \langle u, Y^{(m)} \rangle \mathbf{1}\{Y^{(m)}\notin\gamma_m A\} + \E \beta_m \langle u, \xi \rangle \\
& = \E\alpha_m \E \langle u, Y^{(m)} \rangle -\E\alpha_m \E \langle u, Y^{(m)} \rangle \mathbf{1}\{Y^{(m)}\in \gamma_m A\}
\end{align*}
so that
$$
\E \langle u, Y^{(m)} \rangle = -\frac{\E\alpha_m}{1-\E\alpha_m} \E \langle u,Y^{(m)} \rangle \mathbf{1}\{Y^{(m)}\in\gamma_m A\}.
$$
Together with $\E Y^{(m)}\to\mu$ as $m\to\infty$ and \eqref{alphabeta}, we deduce that
\begin{equation}\label{eqn:Limit_E_Inner_Product}
\lim_{m\to\infty} m \E \langle u, Y^{(m)}\rangle \mathbf{1}\{Y^{(m)}\in \gamma_m A\} = -a \langle u,\mu\rangle.
\end{equation}
Because of the inequality
$$
\big| e^{\mathbf{i} s_1} - \mathbf{i} s_1 - (e^{\mathbf{i} s_2} -\mathbf{i} s_2) \big| = \big| e^{\mathbf{i} s_1} - 1 -\mathbf{i} s_1 - (e^{\mathbf{i} s_2} - 1 -\mathbf{i} s_2) \big| \leq |s_1-s_2|^2, \quad s_1,s_2\in\R,
$$
we obtain
\begin{align*}
& \big| \E \mathbf{1}\{Y^{(m)}\in\gamma_mA\} (e^{\mathbf{i} \alpha_m \langle u,Y^{(m)}\rangle} - e^{\mathbf{i}  \langle u,Y^{(m)}\rangle}) \big| \\
& \leq \E \mathbf{1}\{Y^{(m)}\in\gamma_mA\} (1-\alpha_m)^2 \langle u,Y^{(m)}\rangle^2 + \big| \E \mathbf{1}\{Y^{(m)}\in\gamma_mA\} (1-\alpha_m) \langle u,Y^{(m)}\rangle \big| \\
& \leq \E(1-\alpha_m)^2 \E \mathbf{1}\{Y^{(m)}\in\gamma_mA\}  \langle u,Y^{(m)}\rangle^2 + \big| 1-\E \alpha_m \big| \big| \E \mathbf{1}\{Y^{(m)}\in\gamma_mA\}  \langle u,Y^{(m)}\rangle \big|.
\end{align*}
From \eqref{eqn:Limit_E_Inner_Product} and \eqref{alphabeta} it follows that
$$
\lim_{m\to\infty} m \big| 1-\E \alpha_m \big| \big| \E \mathbf{1}\{Y^{(m)}\in\gamma_mA\}  \langle u,Y^{(m)}\rangle \big|=0,
$$
while the assumption 
\begin{equation}\label{eqn:Limitw}
\lim_{ m\to\infty } m \E \langle u, Y^{(m)} \rangle^2 \mathbf{1}\{Y^{(m)}\in\gamma_m A\} = w,
\end{equation}
and
$$
\lim_{m\to\infty} m \E(\alpha_m-1)^2= \lim_{m\to\infty} m (\E\alpha_m^2-1)-2m(\E\alpha_m-1)=-2a+2a=0,
$$
which is a consequence of \eqref{alphabeta}, lead to 
$$
\lim_{m\to\infty} m \E(1-\alpha_m)^2 \E \mathbf{1}\{Y^{(m)}\in\gamma_mA\}  \langle u,Y^{(m)}\rangle^2 =0.
$$
So we have shown that
$$
\lim_{m\to\infty} m \big| \E \mathbf{1}\{Y^{(m)}\in\gamma_mA\} (e^{\mathbf{i} \alpha_m \langle u,Y^{(m)}\rangle} - e^{\mathbf{i}  \langle u,Y^{(m)}\rangle}) \big| = 0.
$$
This fact, together with \eqref{eqn:Limit_E_Inner_Product} and \eqref{eqn:Limitw}, implies the statement of the lemma if we show that the quantity  
$$
R_m:= m \E \mathbf{1}\{Y^{(m)}\in \gamma_m A\} \big(1 + \mathbf{i}\langle u,Y^{(m)}\rangle - \frac{\langle u,Y^{(m)}\rangle^2}{2} - e^{\mathbf{i}\langle u, Y^{(m)}\rangle}\big), \quad m\in\N,
$$
tends to zero as $m\to\infty$. In the following we consider the cases $\gamma_m\to 0$ and $\gamma_m\to \nu\in(0,\infty]$ as $n\to\infty$ separately. Since the same arguments may be applied to subsequences, this is enough for the proof.

We start with $\gamma_m\to 0$ as $m\to\infty$. The inequality
\begin{equation}\label{eqn:ComplexExponential1}
\big|1+\mathbf{i}s-\frac{s^2}{2}-e^{\mathbf{i}s}\big| \leq \frac{|s|^3}{3}, \quad s\in\R,
\end{equation}
and \eqref{eqn:BoundInnerProductA} lead to
$$
|R_m| \leq \frac{m}{3} \E |\langle u,Y^{(m)}\rangle|^3 \mathbf{1}\{Y^{(m)}\in \gamma_m A\} \leq \frac{\|u\| \overline{r}}{3} \gamma_m m \E \langle u,Y^{(m)}\rangle^2 \mathbf{1}\{Y^{(m)}\in \gamma_m A\}.
$$
Now $\gamma_m\to0$ as $m\to\infty$ and \eqref{eqn:Limitw} yield that $R_m\to0$ as $m\to\infty$.

Next we assume that $\gamma_m\to \nu\in(0,\infty]$ as $m\to\infty$. It follows from Proposition \ref{prop11} (ii) that
\begin{align*}
R_m & = \frac{m}{\E\tau^{(m)}} \E \sum_{j=1}^\infty \mathbf{1}\{\tau^{(m)}=j\} \bigg(1 + \mathbf{i}\langle u, X_j^{(m)}\rangle - \frac{\langle u, X_j^{(m)}\rangle^2}{2} - e^{\mathbf{i}\langle u,X_j^{(m)}\rangle}\bigg) \\
& = \frac{m}{\E\tau^{(m)}} \E \mathbf{1}\{\beta_m \xi \in \gamma_m A\}  \bigg(1+ \mathbf{i} \beta_m \langle u, \xi \rangle - \frac{\beta_m^2 \langle u, \xi \rangle^2}{2} - e^{\mathbf{i}\beta_m\langle u, \xi\rangle }\bigg) \\
& \quad + \frac{m}{\E\tau^{(m)}} \E \sum_{j=2}^\infty \mathbf{1}\{\tau^{(m)}=j\} \mathbf{1}\{X_j^{(m)}\in\gamma_m A\} \bigg(1 + \mathbf{i}\langle u, X_j^{(m)}\rangle - \frac{\langle u, X_j^{(m)}\rangle^2}{2} - e^{\mathbf{i}\langle u,X_j^{(m)}\rangle}\bigg)\\
& =: R_{1,m}+R_{2,m}. 
\end{align*}
Using the inequality
\begin{equation}\label{eqn:ComplexExponential2}
|1+\mathbf{i}s-e^{\mathbf{i}s}|\leq s^2, \quad s\in\R,
\end{equation}
and \eqref{eqn:BoundInnerProductA}, we obtain
\begin{align*}
|R_{2,m}| & \leq  \frac{3}{2} \frac{m}{\E\tau^{(m)}} \E \sum_{j=2}^\infty \mathbf{1}\{\tau^{(m)}=j\} \mathbf{1}\{X_j^{(m)}\in\gamma_m A\} \langle u, X_j^{(m)}\rangle^2 \\
& \leq \frac{3m}{2\E\tau^{(m)}} \|u\|^2 \overline{r}^2 \mathbb{P}(\tau^{(m)}\geq 2) \gamma_m^2 = \frac{3m \|u\|^2 \overline{r}^2}{2\E\tau^{(m)}} \gamma_m^2 \mathbb{P}(\beta_m\xi\notin \gamma_m A) \\
& \leq \frac{3m \|u\|^2 \overline{r}^2}{2\E\tau^{(m)}} \gamma_m^2 \mathbb{P}(\beta_m\|\xi\|\geq \gamma_m \underline{r}) \\
& \leq \frac{3m \|u\|^2 \overline{r}^2}{2\E\tau^{(m)} \underline{r}^2} \E \beta_m^2\|\xi\|^2 \mathbf{1}\{\beta_m\|\xi\|\geq \gamma_m\underline{r}\} = \frac{3\|u\|^2 \overline{r}^2}{2\E\tau^{(m)} \underline{r}^2} m\beta_m^2 \E \|\xi\|^2 \mathbf{1}\{\|\xi\|\geq \gamma_m\underline{r}/\beta_m\}.
\end{align*}
Since, by \eqref{alphabeta} and the assumption on $\{\gamma_m\}_{m\in\N}$, $m\beta_m^2\to 1$ and $\gamma_m/\beta_m\to\infty$ as $m\to\infty$, we see that $R_{2,m}\to0$ as $m\to\infty$. By \eqref{eqn:ComplexExponential1} and \eqref{eqn:ComplexExponential2}, we have
$$
m \mathbf{1}\{\beta_m \xi \in \gamma_m A\} \bigg|1+ \mathbf{i} \beta_m \langle u, \xi \rangle - \frac{\beta_m^2 \langle u, \xi \rangle^2}{2} - e^{\mathbf{i}\beta_m\langle u, \xi\rangle }\bigg| \leq \frac{m\beta_m^3}{3} |\langle u, \xi \rangle|^3
$$
and
$$
m \mathbf{1}\{\beta_m \xi \in \gamma_m A\} \bigg|1+ \mathbf{i} \beta_m \langle u, \xi \rangle - \frac{\beta_m^2 \langle u, \xi \rangle^2}{2}  - e^{\mathbf{i}\beta_m\langle u, \xi\rangle }\bigg| \leq \frac{3}{2} m\beta_m^2 \langle u, \xi \rangle^2
$$  
almost surely. Because of $m\beta_m^3\to 0$ and $m\beta_m^2\to1$ as $m\to\infty$, the dominated convergence theorem implies that $R_{1,m}\to0$ as $m\to\infty$. 
\end{proof}

\subsection{Proof of \eqref{eqn:LimitMixing}-\eqref{eqn:Formula_phi} and 
of \eqref{eqn:vProjection}-\eqref{defB}} 
We know from Proposition \ref{prop:2ndMoment} that $M_2:=\max_{m\in\N}\E \|Y^{(m)}\|^2<\infty$. It follows from the Markov inequality that, for all $z\in(0,\infty)$,
$$
\max_{m\in\N} \mathbb{P}(\|Y^{(m)}\|>z)\leq \frac{M_2}{z^2}.
$$
Thus, the sequence of distributions of $\{Y^{(m)}\}_{m\in\N_0}$ is tight and has a subsequence that converges in distribution. To keep our notation simple we assume that this sequence is the original sequence $\{Y^{(m)}\}_{m\in\N}$ and denote the limiting random vector by $Y$. For $m\in\N$ let $\varphi_{m}$ be the characteristic function of $Y^{(m)}$, while $\varphi$ denotes the characteristic function of $Y$. The weak convergence $Y^{(m)}\overset{d}{\longrightarrow}Y$ as $m\to\infty$ and the uniform boundedness of the second moments imply that
\begin{equation}\label{eqn:ConvergencePhi}
\lim_{m\to\infty} \varphi_m(u)=\varphi(u) \quad \text{and} \quad \lim_{m\to\infty} \varphi_m'(u)=\varphi'(u)
\end{equation}
for $u\in\R^d$. Let $\varphi_\xi$ stand for the characteristic function of $\xi$.

Let $m\in\N$ and $u\in\R^d$ be fixed. We obtain from \eqref{eqn:DefinitionYm} that
\begin{align*}
\varphi_m(u) &  
= \E \mathbf{1}\{Y^{(m)}\notin \gamma_m A\} e^{\mathbf{i}\langle u, \alpha_m Y^{(m)}+\beta_m \xi\rangle} + \E \mathbf{1}\{Y^{(m)}\in \gamma_m A\} e^{\mathbf{i}\beta_m \langle u,\xi\rangle} \\
& = \varphi_\xi(\beta_m u) \big( \E\varphi_m(\alpha_m u) + \E \mathbf{1}\{Y^{(m)}\in \gamma_m A\} (1-e^{\mathbf{i}\alpha_m \langle u,Y^{(m)}\rangle}) \big) \\
& = \varphi_\xi(\beta_m u) \big( \varphi_m(u) + \E\varphi_m(\alpha_m u) - \varphi_m(u) + \E \mathbf{1}\{Y^{(m)}\in \gamma_m A\} (1-e^{\mathbf{i}\alpha_m \langle u,Y^{(m)}\rangle}) \big).
\end{align*}
This can be rewritten as
\begin{equation}\label{eqn:Formula_phi_m}
\frac{1-\varphi_\xi(\beta_mu)}{\varphi_\xi(\beta_mu)} \varphi_m(u) = \E\varphi_m(\alpha_m u) - \varphi_m(u) + \E \mathbf{1}\{Y^{(m)}\in\gamma_m A\} (1-e^{\mathbf{i}\alpha_m \langle u, Y^{(m)}\rangle}).
\end{equation}
Because of $\E\xi=0$ and $\E\|\xi\|^2<\infty$, we have that
$$
\varphi_\xi(s u)=1-\frac{\E\langle u,\xi\rangle^2}{2}s^2 + r(s,u), \quad s\in\R, \quad u\in\R^d,
$$
with $r(s,u)/s^2\to0$ as $s\to0$. Together with \eqref{alphabeta} this yields
\begin{equation}\label{eqn:Quotient_phi_xi}
\lim_{m\to\infty} \frac{m(1-\varphi_\xi(\beta_m u))}{\varphi_\xi(\beta_m u)} = \lim_{m\to\infty} \frac{m\beta_m^2 \E\langle u,\xi\rangle^2}{2} = \frac{\E\langle u,\xi\rangle^2}{2}.
\end{equation}
The Taylor expansion leads to
$$
\E\varphi_m(\alpha_m u) - \varphi_m(u)  = \langle u,\varphi_m'(u)\rangle (\E\alpha_m-1) +\psi_{0,m} \E(\alpha_m-1)^2
$$
with $\psi_{0,m}\in\mathbb{C}$ satisfying $|\psi_{0,m}|\leq \E \langle u, Y^{(m)} \rangle^2\leq \|u\|^2 \E \|Y^{(m)}\|^2 \leq M_2 \|u\|^2$. It follows from \eqref{alphabeta} that 
$$
\lim_{m\to\infty} m (\E\alpha_m-1) = -a
$$
and
$$
\lim_{m\to\infty} m \E(\alpha_m-1)^2= \lim_{m\to\infty} m (\E\alpha_m^2-1)-2m(\E\alpha_m-1)=-2a+2a=0.
$$
Together with \eqref{eqn:ConvergencePhi}, we obtain
\begin{equation}\label{eqn:Diff_phi_m}
\lim_{m\to\infty} m (\E\varphi_m(\alpha_m u) - \varphi_m(u))= \lim_{m\to\infty} -a \langle u,\varphi_m'(u)\rangle=-a \langle u,\varphi'(u)\rangle.
\end{equation}
Lemma \ref{lemm:2ndMoment} and Lemma \ref{lem:LimInd} imply that
\begin{equation}\label{eqn:LimitRemainder}
\lim_{m\to\infty} m \E \mathbf{1}\{Y^{(m)}\in\gamma_m A\} (1-e^{\mathbf{i}\alpha_m \langle u, Y^{(m)}\rangle})=  \mathbf{i} a \langle u,\mu \rangle + \frac{1-p}{2} \E\langle u,\xi\rangle^2,
\end{equation}
where $p$ is defined in \eqref{eqn:Definition_p}. Multiplying \eqref{eqn:Formula_phi_m} with $m$, letting $m\to\infty$, and combining this with \eqref{eqn:ConvergencePhi}, \eqref{eqn:Quotient_phi_xi}, \eqref{eqn:Diff_phi_m}, and \eqref{eqn:LimitRemainder} lead to
$$
\frac{\E\langle u,\xi\rangle^2}{2} \varphi(u)= -a \langle u,\varphi'(u) \rangle + \mathbf{i} a\langle u,\mu\rangle + \frac{1-p}{2}\E\langle u,\xi \rangle^2.
$$
We use the notation 
$$
g(u)=\frac{\sqrt{\E\langle u,\xi\rangle^2}}{\sqrt{2a}}=\frac{\sqrt{u^T\Sigma u}}{\sqrt{2a}}
$$ for $u\in\R^d$ in the sequel. We have shown that $\varphi$ is a solution of the partial differential equation
\begin{equation}\label{eqn:PDE}
\langle u,\varphi'(u)\rangle= - g(u)^2 \varphi(u) + \mathbf{i} \langle u,\mu\rangle + (1-p) g(u)^2, \quad u\in \R^d,
\end{equation}
with $\varphi(0)=1$. Next we show that \eqref{eqn:PDE} has at most one solution. Assume that $\varphi_1$ and $\varphi_2$ are solutions. Then the function $\overline{\varphi}:=\varphi_1-\varphi_2$ satisfies
\begin{equation}\label{eqn:PDEUnique}
\langle u,\overline{\varphi}'(u)\rangle= - g(u)^2 \overline{\varphi}(u), \quad u\in \R^d,
\end{equation}
with $\overline{\varphi}(0)=0$. For a fixed $v\in\R^d$ with $v\neq 0$ we define $\psi(t)=\overline{\varphi}(tv)$, $t\geq 0$. It follows from \eqref{eqn:PDEUnique} and the definition of $g$ that
$$
\psi'(t)= -g(v)^2 \psi(t), \quad t>0,
$$
with $\psi(0)=0$. By the Picard-Lindel\"{o}f theorem this initial value problem has a unique solution, which yields $\psi\equiv 0$. Since we can apply this argument for all choices of $v\neq 0$, we obtain $\varphi_1=\varphi_2$.

In the following, we construct a solution of \eqref{eqn:PDE}. Let $N\sim N(0,1)$ be a standard normal random variable and let $\varphi_{N}$ and $\varphi_{|N|}$ be the characteristic functions of $N$ and $|N|$, respectively, i.e.,
\begin{equation}\label{eqn:phiN_and_phi|N|}
\varphi_N(t)=\exp\bigg(-\frac{t^2}{2} \bigg) \quad \text{and} \quad \varphi_{|N|}(t)=\bigg( 1+\mathbf{i}\frac{\sqrt{2}}{\sqrt{\pi}} \int_0^t\exp\bigg( \frac{x^2}{2}\bigg) \, dx\bigg)\exp\bigg(-\frac{t^2}{2} \bigg), \quad t\in\R.
\end{equation}
Straightforward computations show that 
\begin{equation}\label{eqn:Derivatives}
\varphi'_{N}(t)= -t \varphi_{N}(t), \quad t\in\R, \quad \text{and} \quad \varphi'_{|N|}(t)= \mathbf{i} \frac{\sqrt{2}}{\sqrt{\pi}} - t \varphi_{|N|}(t), \quad t\in\R.
\end{equation}
Moreover, define $\varphi_0(u)=1$, $u\in\R^d$, and let
\begin{align}\label{phiN}
\widetilde{\varphi}(u) = (1-p) \varphi_0(u) + \bigg( p - \frac{\langle u, \mu\rangle}{\sqrt{\frac{2}{\pi}}g(u)} \bigg)\varphi_N(g(u))+ \frac{\langle u, \mu\rangle}{\sqrt{\frac{2}{\pi}}g(u)} \varphi_{|N|}(g(u)).
\end{align}
For a function $f:\R^m\to\mathbb{C}$ and $j\in\{1,\hdots,d\}$ we denote by $D_jf$ its $j$-th partial derivative. It follows from \eqref{eqn:Derivatives} that, for $j\in\{1,\hdots,d\}$,
\begin{align*}
D_j\widetilde{\varphi}(u) & = -\bigg( p - \frac{\langle u, \mu\rangle}{\sqrt{\frac{2}{\pi}}g(u)} \bigg)\varphi_N(g(u)) g(u) D_jg(u) \\
& \quad - \bigg(\frac{\mu_j}{\sqrt{\frac{2}{\pi}}g(u)}-\frac{\langle u, \mu\rangle}{\sqrt{\frac{2}{\pi}}g(u)^2} D_jg(u)\bigg) \varphi_N(g(u)) \\
& \quad + \frac{\langle u, \mu\rangle}{\sqrt{\frac{2}{\pi}}g(u)} \bigg(\mathbf{i}\frac{\sqrt{2}}{\sqrt{\pi}} - g(u) \varphi_{|N|}(g(u)) \bigg) D_jg(u)\\
& \quad + \bigg(\frac{\mu_j}{\sqrt{\frac{2}{\pi}}g(u)}-\frac{\langle u, \mu\rangle}{\sqrt{\frac{2}{\pi}}g(u)^2} D_jg(u)\bigg) \varphi_{|N|}(g(u)).
\end{align*}
This implies that
\begin{align*}
\langle u,\widetilde{\varphi}'(u) \rangle & = - \langle u,g'(u)\rangle g(u) \bigg( \bigg( p - \frac{\langle u, \mu\rangle}{\sqrt{\frac{2}{\pi}}g(u)} \bigg)\varphi_N(g(u)) + \frac{\langle u, \mu\rangle}{\sqrt{\frac{2}{\pi}}g(u)} \varphi_{|N|}(g(u)) \bigg)\\
& \quad + \mathbf{i}\frac{\langle u,\mu\rangle \langle u, g'(u)\rangle}{g(u)} - \bigg(\frac{\langle u,\mu\rangle}{\sqrt{\frac{2}{\pi}}g(u)}-\frac{\langle u, \mu\rangle}{\sqrt{\frac{2}{\pi}}g(u)^2} \langle u,g'(u)\rangle \bigg) (\varphi_N(g(u)) - \varphi_{|N|}(g(u))).
\end{align*}
Since $g(tu)=tg(u)$ for $t>0$, taking the derivative with respect to $t$ for $t=1$ yields $\langle u,g'(u)\rangle = g(u)$. Thus the previous equation simplifies to
\begin{align*}
\langle u,\widetilde{\varphi}'(u) \rangle & = - g(u)^2 \bigg( \bigg( p - \frac{\langle u, \mu\rangle}{\sqrt{\frac{2}{\pi}}g(u)} \bigg)\varphi_N(g(u)) + \frac{\langle u, \mu\rangle}{\sqrt{\frac{2}{\pi}}g(u)} \varphi_{|N|}(g(u)) \bigg) + \mathbf{i} \langle u,\mu\rangle \\
& \quad - (1-p) g(u)^2 \varphi_0(u) + (1-p) g(u)^2 \\
& = -g(u)^2\widetilde{\varphi}(u) + \mathbf{i} \langle u,\mu\rangle + (1-p) g(u)^2. 
\end{align*}
Since $\widetilde{\varphi}(0)=1$ and $\eqref{eqn:PDE}$ has at most one solution, this proves that $\varphi=\widetilde{\varphi}$. It follows from \eqref{eqn:phiN_and_phi|N|} that
$$
\varphi_{|N|}(t)= \bigg( 1+\mathbf{i}\frac{\sqrt{2}}{\sqrt{\pi}} \int_0^t\exp\bigg( \frac{x^2}{2}\bigg) \, dx\bigg)\varphi_{N}(t), \quad t\in\R.
$$
This implies that
\begin{align*}
\widetilde{\varphi}(u)& = 1-p + \bigg( p - \frac{\langle u, \mu\rangle}{\sqrt{\frac{2}{\pi}}g(u)} \bigg)\varphi_N(g(u))\\
& \quad + \frac{\langle u, \mu\rangle}{\sqrt{\frac{2}{\pi}}g(u)} \bigg( 1+\mathbf{i}\frac{\sqrt{2}}{\sqrt{\pi}} \int_0^{g(u)} \exp\bigg( \frac{x^2}{2}\bigg) \, dx\bigg)\varphi_{N}(g(u)) \\
& = 1-p + \bigg( p +\mathbf{i} \frac{\sqrt{2a}\langle u,\mu\rangle}{\sqrt{u^T\Sigma u}} \int_0^{g(u)} \exp\bigg( \frac{x^2}{2}\bigg) \, dx\bigg) \exp\bigg(-\frac{u^T\Sigma u}{4a}\bigg)
\end{align*}
and, thus, proves \eqref{eqn:Formula_phi}. Using the obvious identity $\varphi_N(t)=\frac{1}{2} (\varphi_{|N|}(t)+\varphi_{-|N|}(t))$, $t\in\R$ in \eqref{phiN}, we get 
$$
\widetilde{\varphi}(u)= 1-p + \bigg(\frac{p}{2} - \frac{\langle u, \mu\rangle}{2\sqrt{\frac{2}{\pi}}g(u)} \bigg) \varphi_{-|N|}(g(u)) + \bigg(\frac{p}{2} + \frac{\langle u, \mu\rangle}{2\sqrt{\frac{2}{\pi}}g(u)} \bigg)\varphi_{|N|}(g(u)), \quad u\in\R^d.
$$
For a fixed $v\in\R^d$, letting $u=tv$ shows that $\langle v,Y \rangle$ has the characteristic function
$$
\varphi_{\langle v,Y \rangle}(t) = 1-p + \bigg(\frac{p}{2} - \frac{\langle v, \mu\rangle}{2\sqrt{\frac{2}{\pi}}g(v)} \bigg) \varphi_{-g(v)|N|}(t) + \bigg(\frac{p}{2} + \frac{\langle v, \mu\rangle}{2\sqrt{\frac{2}{\pi}}g(v)} \bigg)\varphi_{g(v)|N|}(t), \quad t\in\R.
$$
This identity is of the form
$$
\varphi_{\langle v,Y \rangle}(t) = 1-p + c_{-} \varphi_{-g(v)|N|}(t) + c_+\varphi_{g(v)|N|}(t).
$$
We have that $c_-,c_+\geq 0$. To prove this, we assume that $c_-<0$ (the case $c_+<0$ goes analogously). Since $1-p+c_{-}+c_{+}=1$ and $1-p\leq 1$, we have that $c_+> 0$. The above equation can be rewritten as
$$
\frac{1}{1+|c_-|}(\varphi_{\langle v,Y \rangle}(t) + |c_-| \varphi_{-g(v)|N|}(t)) = \frac{1}{1+|c_-|} \bigg( 1-p + c_+\varphi_{g(v)|N|}(t) \bigg).
$$
Here, the left- and the right-hand sides are characteristic functions of mixtures of random variables. Since the distribution belonging to the left-hand side allows negative values and the one belonging to the right-hand side is non-negative, this is a contradiction. So $\langle v,Y\rangle$ is a mixture of the random variables $0$, $-g(v)|N|$, and $g(v) |N|$ and we have shown \eqref{eqn:vProjection}.
\qed

\subsection{Proof of \eqref{densZ}-\eqref{denseve}} 

Note that $\varphi_Z$ can be written as
$$
\varphi_Z(u) = \bigg( \frac{1}{2} - \frac{\sqrt{\pi a} \langle \mu, u \rangle}{2p\sqrt{u^T\Sigma u}} \bigg) \varphi_{-\frac{|N|}{\sqrt{2a}}}(\sqrt{u^T\Sigma u}) + \bigg( \frac{1}{2} + \frac{\sqrt{\pi a} \langle \mu, u \rangle}{2p\sqrt{u^T\Sigma u}} \bigg) \varphi_{\frac{|N|}{\sqrt{2a}}}(\sqrt{u^T\Sigma u}), \quad u\in\R^d.
$$
For $u\in\R^d$ this implies that $\langle u, Z\rangle$ has the density
$$
f_{u,\mu,\Sigma}(s) = \bigg( \frac{1}{2} + \frac{s}{|s|} \frac{\sqrt{\pi a} \langle \mu, u \rangle}{2p\sqrt{u^T\Sigma u}} \bigg) \frac{\sqrt{2}}{\sqrt{\pi}} \frac{\sqrt{2a}}{\sqrt{u^T\Sigma u}} \exp\bigg(-\frac{a s^2}{u^T\Sigma u}\bigg), \quad s\in\R.
$$
We first prove the statements for $\Sigma=I_d$, where $I_d$ denotes the $d$-dimensional identity matrix, and for odd $d$. For $d=1$, $f_{1,\mu,I_1}$ simplifies to $f_Z$. Thus, it is sufficient to assume that $d=2k+3$ with $k\in\N_0$ in the sequel. Let $\widetilde{Z}$ be a random vector in $\R^d$ with density
$$
g_{\mu}(x)= \frac{\sqrt{a}^d}{\sqrt{\pi}^d } e^{-a\|x\|^2} +   \frac{(-2)^{k+1} a^{k+\frac{5}{2}}}{p\kappa_{2k+2}(2k+2)!! } \langle \mu, x\rangle h^{(k+1)}(a\|x\|^2), \quad x\in\R^d.
$$
For a fixed $u\in\R^d$ we compute the density $g_{u,\mu}$ of $\langle u,\widetilde{Z}\rangle$ in the following. For $s\in\R$ we obtain that
\begin{align*}
g_{u,\mu}(s) & = \lim_{\varepsilon\to 0} \frac{\mathbb{P}(\langle u,\widetilde{Z}\rangle \leq s+\varepsilon) - \mathbb{P}(\langle u,\widetilde{Z}\rangle \leq s)}{\varepsilon} \\
& =   \lim_{\varepsilon\to 0} \frac{1}{\varepsilon} \int_{\{z\in\R^d: s\leq \langle u,z \rangle \leq s+\varepsilon\}} g_{\mu}(y) \, dy = \frac{1}{\|u\|} \int_{H_{u,s}} g_{\mu}(y) \, dy
\end{align*}
with $H_{u,s}:=\{z\in\R^d: \langle u,z \rangle=s\}$. This can be rewritten as
$$
g_{u,\mu}(s)= \frac{1}{\|u\|} \int_{H_{u,0}} g_{\mu}\bigg(\frac{s}{\|u\|^2}u+y\bigg) \, dy.
$$
Since $\langle u, y\rangle=0$ for all $y\in H_{u,0}$, we obtain that
\begin{align*}
g_{u,\mu}(s) & = \frac{1}{\|u\|} \int_{H_{u,0}} \frac{\sqrt{a}^d}{\sqrt{\pi}^d} \exp\bigg(-a \bigg(\frac{s^2}{\|u\|^2}+\|y\|^2\bigg)\bigg) \\
& \quad \quad \quad \quad \quad \quad + \frac{(-2)^{k+1} a^{k+\frac{5}{2}}}{p\kappa_{2k+2}(2k+2)!!} \langle \mu, \frac{s}{\|u\|^2}u+y\rangle h^{(k+1)}\bigg(\frac{as^2}{\|u\|^2}+a\|y\|^2\bigg) \, dy \\
& = \frac{1}{\|u\|} \int_{H_{u,0}} \frac{\sqrt{a}^d}{\sqrt{\pi}^d} \exp\bigg(-\bigg(\frac{as^2}{\|u\|^2}+a\|y\|^2\bigg)\bigg) \, dy \\
& \quad + \frac{(-2)^{k+1} a^{k+\frac{5}{2}}}{p\kappa_{2k+2}(2k+2)!!} \frac{\langle \mu,u\rangle }{\|u\|^3} s \int_{H_{u,0}} h^{(k+1)}\bigg(\frac{as^2}{\|u\|^2}+a\|y\|^2\bigg) \, dy \\
& \quad + \frac{(-2)^{k+1} a^{k+\frac{5}{2}}}{p\kappa_{2k+2}(2k+2)!!} \frac{1}{\|u\|} \int_{H_{u,0}} \langle \mu, y\rangle h^{(k+1)}\bigg(\frac{as^2}{\|u\|^2}+a\|y\|^2\bigg) \, dy \\
& =: G_1+G_2+G_3.
\end{align*}
A straightforward computation shows that
$$
G_1 = \frac{\sqrt{a}}{\sqrt{\pi}\|u\|} \exp\bigg(-\frac{as^2}{\|u\|^2}\bigg).
$$
From the fact that the integrand is an odd function, it follows that $G_3=0$. Using polar coordinates, we obtain
\begin{align*}
\int_{H_{u,0}} h^{(k+1)}\bigg(\frac{as^2}{\|u\|^2}+a\|y\|^2\bigg) \, dy & = \int_{\R^{d-1}} h^{(k+1)}\bigg(\frac{as^2}{\|u\|^2}+a\|y\|^2\bigg) \, dy \\
& = (d-1)\kappa_{d-1} \int_0^\infty h^{(k+1)}\bigg(\frac{as^2}{\|u\|^2}+ar^2\bigg) r^{d-2} \, dr \\
& = (2k+2)\kappa_{2k+2} \int_0^\infty h^{(k+1)}\bigg(\frac{as^2}{\|u\|^2}+ar^2\bigg) r^{2k+1} \, dr.
\end{align*}
Iterating integration by parts $k$ times leads to
\begin{align*}
\int_0^\infty h^{(k+1)}\bigg(\frac{as^2}{\|u\|^2}+ar^2\bigg) r^{2k+1} \, dr & = \frac{(-1)^k (2k)!!}{(2a)^{k}} \int_0^\infty h'\bigg(\frac{as^2}{\|u\|^2}+ar^2\bigg) r \, dr \\
& = \frac{(-1)^{k+1} (2k)!!}{(2a)^{k+1}} h\bigg(\frac{as^2}{\|u\|^2}\bigg).
\end{align*}
This shows that
\begin{align*}
G_2 & = a^{\frac{3}{2}} \frac{\langle \mu, u \rangle}{p\|u\|^3} s  h\bigg(\frac{as^2}{\|u\|^2}\bigg) = a^{\frac{3}{2}} \frac{\langle \mu,u \rangle}{p\|u\|^3} s \frac{\exp\big(-\frac{as^2}{\|u\|^2}\big)}{\sqrt{a \frac{s^2}{\|u\|^2}}} = a \frac{s}{p|s|} \frac{\langle \mu,u \rangle}{\|u\|^2} \exp\bigg(-\frac{as^2}{\|u\|^2}\bigg) \\
& = \frac{s}{|s|} \frac{\sqrt{\pi a}}{2p} \frac{\langle \mu,u \rangle}{\|u\|} \frac{\sqrt{2}}{\sqrt{\pi}} \frac{\sqrt{2a}}{\|u\|} \exp\bigg(-\frac{as^2}{\|u\|^2}\bigg),
\end{align*}
whence
$$
g_{u,\mu}(s) = \bigg( \frac{1}{2} + \frac{s}{|s|} \frac{\sqrt{\pi a}}{2p} \frac{\langle \mu,u \rangle}{\|u\|} \bigg) \frac{\sqrt{2}}{\sqrt{\pi}} \frac{\sqrt{2a}}{\|u\|} \exp\bigg(-\frac{as^2}{\|u\|^2}\bigg).
$$
This means that $g_{u,\mu}=f_{u,\mu,I_d}$. Thus, $Z$ has the density $g_\mu$, which proves (ii) for $\Sigma=I_d$ and $d$ odd.

Next we consider the case that $d$ is even and $\Sigma=I_d$. In order to emphasize the dependence on $\mu$ and $\Sigma$, we write $Z=Z_{\mu,\Sigma}$ and denote the density of $Z$ by $f_{\mu,\Sigma}$ and the characteristic function by $\varphi_{\mu,\Sigma}$. Note that $Z_{\mu,\Sigma}$ follows the same distribution as the first $d$ components of $Z_{(\mu,0),I_{d+1}}$, where $(\mu,0)$ is the vector in $\R^{d+1}$ such that the first $d$ components coincide with those of $\mu$ and the last component is zero. Since $Z_{(\mu,0),I_{d+1}}$ has the density $f_{(\mu,0),I_{d+1}}$, we have
$$
f_{\mu,I_d}(x)=\int_{-\infty}^{\infty} f_{(\mu,0),I_{d+1}}((x,z)) \, dz,
$$ 
which proves (ii) for $\Sigma=I_d$ and $d$ even.

Finally, we consider the case of a general positive definite covariance matrix $\Sigma\in\R^{d\times d}$. Let $\Sigma^{-\frac{1}{2}}$ be the unique positive definite matrix in $\R^{d\times d}$ such that $\Sigma^{-\frac{1}{2}} \Sigma^{-\frac{1}{2}}=\Sigma^{-1}$. For $u\in\R^d$ we have that 
\begin{align*}
\int_{\R^d} e^{\mathbf{i} \langle u,x \rangle} \frac{1}{\sqrt{\det(\Sigma)}} f_{\Sigma^{-\frac{1}{2}}\mu,I_d}(\Sigma^{-\frac{1}{2}}x) \, dx & = \int_{\R^d} e^{\mathbf{i} \langle \Sigma^{\frac{1}{2}}u, y \rangle} f_{\Sigma^{-\frac{1}{2}}\mu,I_d}(y) \, dy \\
& = \varphi_{\Sigma^{-\frac{1}{2}}\mu,I_d}(\Sigma^{\frac{1}{2}}u) = \varphi_{\mu,\Sigma}(u),
\end{align*}
whence $f_{\mu,\Sigma}(\cdot)=\frac{1}{\sqrt{\det(\Sigma)}} f_{\Sigma^{-\frac{1}{2}}\mu,I_d}(\Sigma^{-\frac{1}{2}}\cdot)$.
\qed

\subsection{Proof of Proposition \ref{thm:NoTruncation}} 

It follows from Theorem 1.5 and Theorem 1.6 in \cite{Vervaat1979} that $X^{(m)}$ has the same distribution as
$$
\beta_m \sum_{i=1}^\infty \prod_{j=1}^{i-1} \alpha_{m,j}\xi_i.
$$
This leads to
$$
\E \|X^{(m)}\|^2 = \beta_m^2 \sum_{i=1}^\infty (\E \alpha_m^2)^{i-1} \E \|\xi\|^2 = \frac{\beta_m^2}{1-\E\alpha_m^2} \E\|\xi\|^2,
$$
whence, by \eqref{alphabeta},
$$
\lim_{m\to\infty} \E \|X^{(m)}\|^2 =\frac{1}{2a}.
$$
Thus, $\{X^{(m)}\}_{m\in\N}$ is tight. By following the arguments of step 2 of the proof of Theorem \ref{thm:LimitMixing} for $A=\varnothing$, we obtain that the characteristic function $\varphi$ of the limiting distribution of a weakly convergent subsequence of $\{X^{(m)}\}_{m\in\N}$ satisfies
$$
\frac{\E\langle u,\xi \rangle^2}{2}\varphi(u) = -a \langle u,\varphi'(u) \rangle, \quad u\in\R^d.
$$
The observation that the characteristic function of $N_{\widehat{\Sigma}}$ is the unique solution of this differential equation completes the proof. \qed

\appendix

\section{On the proof of Proposition \ref{prop11}}\label{AppendixA}

To make the paper self-contained, we provide here our proof
of part (i) of the proposition.
Part (ii) may be found, say,  in \cite{Asm} (page 171, Corollary VI.1.4).
Then part (iii) follows, say, from Theorem 1.4 in \cite{AlsIksRoe} (see also further references therein). Indeed, 
the stationary distribution of $\{||Y_t^{(m)}||\}_{t\in\N_0}$ is stochastically dominated by that
of the one-dimensional recursion 
\begin{align*}
\widehat{Y}^{(m)}_{t+1}= \alpha_{m,t+1} \widehat{Y}^{(m)}_{t} \mathbf{1} \{ 
\widehat{Y}^{(m)}_{t} > \gamma_m \underline{r}\} + \beta_m ||\xi_{t+1}||, \quad t=0,1,\ldots,
\end{align*}
which, in turn, is stochastically dominated by the stationary
distribution of the one-dimensional recursion
\begin{align*}
\widehat{X}^{(m)}_{t+1} = \alpha_{m,t+1} \widehat{X}^{(m)}_{t} + \beta_m||\xi_{t+1}||, \quad t=0,1,\ldots,
\end{align*}
and the latter distribution was considered in \cite{AlsIksRoe}. Recall that, for two
one-dimensional random variables $Z_1$ and $Z_2$, we say that the distribution of $Z_1$ is stochastically dominated by the distribution of $Z_2$ if ${\mathbb P} (Z_1>x) \le {\mathbb P} (Z_2>x)$, for all $x\in\R$.

Our proof of  part (i) is based on two lemmas that may be of an independent interest.

\begin{lemm}\label{REC} 
Assume that $\E\xi=0$ and that $\E \xi\xi^T$ exists and is regular. Let $\alpha$ be a positive random variable such that $\mathbb{P}(\alpha\in[d/(d+1),1))>0$. Let $Z_{t+1} = \alpha_{t+1} Z_t + \beta \xi_{t+1}$, $t\in\N_0$ with initial value $Z_0$, where $\{\alpha_t\}_{t\in\N}$ and $\{\xi_t\}_{t\in\N}$ are i.i.d.\ copies of $\alpha$ and $\xi$, respectively, and $\beta >0$. Then, for any $\gamma > 0$ and any $R>\gamma$, there exists an integer $L>0$ such that
\begin{align*}
\inf_{x\in B^d(0,R)} {\mathbb P} (\|Z_t\|\le \gamma \ \mbox{for some} \
t\le L \ | \ Z_0=x ) >0.
\end{align*}
\end{lemm}

The following example illustrates that $\alpha$ may not be too small. It can be extended to higher dimensions.

\noindent {\bf Example.} Let $\alpha$ be a constant in $(0,1/2)$ and let $\beta=1$ and $Z_0=0$. Assume that $\mathbb{P}(\xi=1)=\mathbb{P}(\xi=-1)=\frac{1}{2}$. For $t\in\N$ we have
$$
Z_t=\sum_{i=1}^t \alpha^{t-i} \xi_i \quad \text{and} \quad \big| \sum_{i=1}^{t-1} \alpha^{t-i} \xi_i \big| \leq \sum_{i=1}^{t-1} \alpha^{t-i} \leq \frac{\alpha}{1-\alpha}
$$
so that
$$
|Z_t|\geq |\xi_t| - \frac{\alpha}{1-\alpha} = 1 - \frac{\alpha}{1-\alpha} =\frac{1-2\alpha}{1-\alpha}.
$$
This shows that $\{Z_t\}_{t\in\N}$ cannot hit an interval $[-\gamma,\gamma]$ for $\gamma\in (0, (1-2\alpha)/(1-\alpha))$.

\begin{proof}[Proof of Lemma \ref{REC}]
We can assume $\beta=1$ without loss of generality. The supports of $\alpha$ and $\xi$ are given by
$$
\operatorname{supp} \alpha = \{x\in\R: \mathbb{P}(\alpha\in (x-r,x+r) )>0 \text{ for all } r>0\}
$$
and
$$
\operatorname{supp} \xi = \{x\in\R^d: \mathbb{P}(\xi\in B^d(x,r))>0 \text{ for all } r>0\}.
$$
Let $\operatorname{conv}(M)$ denote the convex hull of a set $M\subseteq \R^d$. Note that $0\in\operatorname{conv}(\operatorname{supp} \xi)$. Otherwise there exists a separating hyperplane between $0$ and $\operatorname{conv}(\operatorname{supp} \xi)$. Since this hyperplane can be chosen through the origin, there exists a $w\in\R^d$ such that $\langle w, x\rangle\geq 0$ for all $x\in\operatorname{supp} \xi$, whence $\langle w,\xi \rangle \geq 0$ a.s. We have that $\mathbb{P}(\langle w,\xi \rangle >0)>0$ since otherwise $\langle w,\xi \rangle = 0$ a.s.\ and $\E \xi\xi^T$ would not be regular. This yields $\mathbb{E}\langle w,\xi \rangle > 0$, which contradicts $\mathbb{E}\xi=0$.

The assumption $\mathbb{P}(\alpha\in[d/(d+1),1))>0$ allows us to choose $s\in \operatorname{supp} \alpha \cap [d/(d+1),1)$.
Because of $0\in\operatorname{conv}(\operatorname{supp} \xi)$, there exist $u_0,\hdots,u_d\in \operatorname{supp} \xi$ such that $\sum_{i=0}^d \lambda_i u_i= 0$ for some $\lambda_0,\hdots,\lambda_d\in[0,1]$ with $\sum_{i=0}^d\lambda_i=1$. We define a sequence $(k_\ell)_{\ell\in\mathbb{N}}$ with values in $\{0,\hdots,d\}$ recursively by
$$
k_{\ell}=\operatorname{argmax}_{i\in\{0,\hdots,d\}} \lambda_i-(1-s)\sum_{v=1}^{\ell-1} \mathbf{1}\{k_v=i\} s^{v-1}, \quad \ell\in\N.
$$
In case that the $\operatorname{argmax}$ is not unique, we choose the smallest solution. We prove now by induction over $\ell$ that
\begin{equation}\label{eqn:Inequality_lambda_i}
\lambda_i-(1-s)\sum_{v=1}^{\ell-1} \mathbf{1}\{k_v=i\} s^{v-1}\geq 0
\end{equation}
for all $i\in\{0,\hdots,d\}$ and $\ell\in\N$. For $\ell=1$ this is obvious. We have
$$
\sum_{i=0}^d \bigg( \lambda_i-(1-s)\sum_{v=1}^{\ell-1} \mathbf{1}\{k_v=i\} s^{v-1} \bigg) = 1 - (1-s)\sum_{v=1}^{\ell-1}s^{v-1} = (1-s)\sum_{v=\ell}^{\infty}s^{v-1}=s^{\ell-1}
$$
so that
$$
\max_{i\in\{0,\hdots,d\}} \lambda_i-(1-s)\sum_{v=1}^{\ell-1} \mathbf{1}\{k_v=i\} s^{v-1} \geq \frac{s^{\ell-1}}{d+1}.
$$
Then $s\in[d/(d+1),1]$ yields
$$
\frac{s^{\ell-1}}{d+1} \geq (1-s)s^{\ell-1},
$$
whence together with the induction assumption
$$
\lambda_i-(1-s)\sum_{v=1}^{\ell} \mathbf{1}\{k_v=i\} s^{v-1}\geq 0
$$
for all $i\in\{0,\hdots,d\}$. Combining \eqref{eqn:Inequality_lambda_i} with
$$
\sum_{i=0}^d \bigg( \lambda_i-(1-s)\sum_{v=1}^{\infty} \mathbf{1}\{k_v=i\} s^{v-1} \bigg)=0
$$
shows that
$$
(1-s)\sum_{\ell=1}^\infty \mathbf{1}\{k_\ell=i\} s^{\ell-1}=\lambda_i
$$
for $i\in\{0,\hdots,d\}$. This implies that
$$
\sum_{\ell=1}^\infty s^{\ell-1} u_{k_\ell} = \sum_{i=0}^d \sum_{\ell=1}^\infty \mathbf{1}\{k_\ell=i\} s^{\ell-1} u_i = \frac{1}{1-s} \sum_{i=0}^d \lambda_i u_i=0.
$$
Since, for $t\in\N$ sufficiently large,
$$
\bigg\|\sum_{\ell=1}^t s^{\ell-1} u_{k_\ell}\bigg\|\leq \frac{\gamma}{3} \quad \text{and} \quad s^t R \leq \frac{\gamma}{3},
$$
there exists an $t_0\in\N$ such that 
$$
\max_{x\in B^d(0,R)}\bigg\|s^{t_0}x + \sum_{\ell=1}^{t_0} s^{\ell-1} u_{k_\ell}\bigg\|\leq \frac{2\gamma}{3}.
$$
We can even choose $\delta>0$ such that
$$
\max_{x\in B^d(0,R)}\bigg\| \prod_{i=1}^{t_0}s_i x + \sum_{\ell=1}^{t_0} \prod_{i=1}^{\ell-1} s_i u_{\ell} \bigg\|\leq \gamma
$$
for all $u_j\in B^d(u_{k_j},\delta)$ and $s_j\in(s-\delta,s+\delta)$, $j\in\{1,\hdots,t_0\}$. So we have shown that, for any $x\in B^d(0,R)$,
$$
\mathbb{P}( \|Z_{t_0}\|\leq \gamma \ | \ Z_0=x) \geq \mathbb{P}(\xi_i\in B^d(u_{k_{t_0+1-i}},\delta) \text{ and } \alpha_i\in(s-\delta,s+\delta), i\in\{1,\hdots,t_0\}).
$$
Because of $u_0,\hdots,u_d\in\operatorname{supp}\xi$ and $s\in\operatorname{supp}\alpha$, the probability on the right-hand side is positive, which completes the proof.
\end{proof}

\begin{lemm}\label{largesmall}
Let $\{X_t\}_{t\in\N_0}$ be a (time-homogeneous) Markov chain taking values in a measurable state space $({\cal X, B_X})$ with initial value $X_0$ and let $\|\cdot \|$ be a non-negative test function. For $R>0$, let
$A_R = \{  x\in {\cal X}: \|x\|\le R\}$ be a ``ball'' of $\|\cdot \|$-radius $R$.
Let $0 < \gamma <R$ be another number.
Assume that\\
(i) the set $A_R$ is {\it geometrically positive recurrent}, i.e.
\begin{align*}
\tau_x (A_R) = \min \{t\ge 1: X_t \in A_R \ | \ X_0=x\}
\end{align*}
is finite a.s. for any $x\in {\cal X}$ and, for some $c_0>0$, 
\begin{align*}
g(c_0) :=\sup_{x\in A_R} \E e^{c_0\tau_x(A_R)} < \infty
\end{align*}
and that\\
(ii) there exists $L\ge 1$ such that
\begin{align*}
\delta = \inf_{x\in A_R } \delta_x >0
\quad
\mbox{where}
\quad
\delta_x = {\mathbb P} (\tau_x(A_{\gamma})\le L) \quad \mbox{and}
\quad A_{\gamma} = \{ x\in {\cal X}: \|x\|\le \gamma\}.
\end{align*}
Then the set $A_{\gamma}$ is geometrically positive recurrent too.
\end{lemm}

\begin{proof}
In order to avoid trivialities, let $\delta <1$.

Since $\tau_x(A_R)$ is a.s. finite for all $x\in \cal{X}$ and $\delta_y \ge \delta >0$ for any $y\in A_R$,
one can easily deduce that $\tau_x(A_{\gamma})$ is a.s. finite for all $x\in \cal{X}$, too.
Now we prove that $\sup_{x_0\in A_{\gamma}} {\mathbb E} e^{c\tau_{x_0}(A_{\gamma})}$ is finite too.

Assume that $X_0=x_0\in A_{\gamma}$ and 
let $V_0=T_0=0$ 
and, for $k=0,1,\ldots,$
\begin{align*}
V_{k+1} = \min \big\{ T_{k}+L, \min \{t>T_{k}: \  \|X_{t}\| \le
 \gamma\} \big\}, \quad 
T_{k+1} = \min \{ t> V_{k+1} : \ \|X_{t}\|\le R\}.
 \end{align*}
Define $\mu = \min \{k\ge 1: \ \|X_{V_k}\| \le \gamma\}$. Clearly, $V_{\mu}\ge \tau_{x_0}(A_{\gamma})$ and 
$T_{\mu -1}\le V_{\mu} < T_{\mu}$ a.s. Further, for any fixed $k$,
\begin{align*}
T_{k+1}-T_k \le \sum_{i=1}^{L+1} \tau_i
\end{align*}
where $\tau_1 = \min \{t>T_k: \|X_t\|\le R\} - T_k$ and, for
$j=1,2,\ldots$, $\tau_{j+1} = \min \{t> T_k+\sum_{\ell=1}^j \tau_\ell: 
\|X_t\|\le R\} - T_k - \sum_{\ell=1}^j \tau_\ell$. By the lemma assumption, for $c\in [0,c_0]$,
\begin{align}\label{supk}
\sup_{x\in A_R}{\mathbb E} \left(e^{c(T_{k+1}-T_k)} \ | \ X_{T_k}=x \right) \le 
g(c)^{L+1}  < \infty
\end{align}
where $g(c)\to 1$ as $c\downarrow 0$.

Now introduce a geometric random variable $\nu$ on the common probability space with $\{X_t\}_{t\in\N_0}$. 
For any $k=0,1,\ldots$, define the $0$-$1$-valued random variable $\theta_{k+1}$ as follows: given
$X_{T_k}=x$, we let $\theta_{k+1}=0$ if $\|X_t\|> \gamma$ 
for all $t=T_k+i$, $i=1,2,\ldots,L$. Otherwise,
we let $\theta_{k+1}=1$ with probability $\delta/\delta_x$ and
$0$ with probability $1-\delta/\delta_x$ independently of everything else. Then we let $\nu = \min \{k: \theta_k=1\}$. Since $\{\theta_k\}_{k\in\N}$ form an i.i.d.\ sequence and
${\mathbb P} (\theta_k=1) = 1-{\mathbb P}(\theta_k=0)=\delta$, the random variable $\nu$ has a geometric distribution with parameter $\delta$. 
Clearly, $\mu \le \nu$ a.s. and $V_{\mu}\le T_{\mu} \le T_{\nu}$ a.s.

For any $c>0$, 
\begin{align*}
{\mathbb E} (e^{cT_{\nu}} \ | \ X_0=x_0)
=
\sum_{k=1}^{\infty} 
{\mathbb E} \left( e^{cT_k}{\mathbf 1} \{\nu=k\} \ | \ X_0=x_0\right). 
\end{align*}
Since, for any $x\in A_R$ and for $k=0,1,\ldots$,  
\begin{align*}
{\mathbb E} \left(e^{c(T_{k+1}-T_{k})} \ | \ 
X_{T_{k}}=x \right) &= {\mathbb P} (\theta_{k+1}=0)
{\mathbb E} \left(e^{c(T_{k+1}-T_{k})} \ | \ 
X_{T_{k}}=x, \theta_{k+1}=0 \right)\\
 & \quad + 
{\mathbb P} (\theta_{k+1}=1)
{\mathbb E} \left(e^{c(T_{k+1}-T_{k})} \ | \ 
X_{T_{k}}=x, \theta_{k+1}=1 \right),
\end{align*}
inequality \eqref{supk} implies that, for any $c\in (0,c_0]$, 
$$
g_0(c)= \sup_{x\in A_R} {\mathbb E} \left(e^{c(T_{k+1}-T_{k})} \ | \ 
X_{T_{k}}=x, \theta_{k+1}=0 \right)
$$
and
$$  
g_1(c)= \sup_{x\in A_R} {\mathbb E} \left(e^{c(T_{k+1}-T_{k})} \ | \ 
X_{T_{k}}=x, \theta_{k+1}=1 \right)
$$
are finite and $g_i(c) \downarrow 1$ as $c\downarrow 0$, for $i=0,1$.

Choose $c\in (0,c_0]$ such that $g_0(c) < (1-\delta)^{-1}$. Then, 
using the strong Markov property,  conditioning and the backward induction argument, one can conclude that 
\begin{align*}
& {\mathbb E} \left( e^{cT_k}{\mathbf 1}\{\nu=k\} \ | \ X_0=x_0\right) \\
&= {\mathbb E} \left(\prod_{i=1}^{k-1} \left(e^{c(T_i-T_{i-1})}{\mathbf 1}\{\theta_i=0\}\right) \cdot
e^{c(T_k-T_{k-1})}{\mathbf 1}(\theta_k=1) \ | \ X_0=x_0\right)\\
&\le (1-\delta)^{k-1}\delta \cdot g_0(c)^{k-1}g_1(c),
\end{align*}
so ${\mathbb E} (e^{cT_{\nu}} \ | \ X_0=x_0) \le \delta g_1(c) (1-(1-\delta)g_0(c))^{-1} \equiv C$ for
all $x_0\in A_{\gamma}$. Therefore, ${\mathbb E} e^{\tau_{x_0}(A_{\gamma})} \le C$ for all $x_0\in A_{\gamma}$, and the result follows.
\end{proof}

Let $\{X_{t,x}^{(m)}\}_{t\in\N_0}$ be an autoregressive sequence $X_{t+1,x}^{(m)} = \alpha_{m,t+1} X_{t,x}^{(m)}+\beta_m \xi_{t+1}$ with initial value $X_{0,x}^{(m)}=x$. To complete the proof of the first statement of Proposition \ref{prop11}, it remains to show that, for some $R>1$, 
the random variables
\begin{align*}
T_x\equiv T_x(R) = \min \{t>0: \ \|X_{t,x}^{(m)}\| \leq R \}
\end{align*} 
(these are the return times to the closure of the ball $B^d(0,R)$) have a uniformly finite exponential moment,
\begin{align}\label{uniRW}
\sup_{x\in A_R} {\mathbb E} e^{cT_x} <\infty, \ \  \mbox{for some}  \ \ c>0.
\end{align}

We prove \eqref{uniRW} now. We may assume, without loss of generality, that $\xi$ is a one-dimensional non-negative random variable and that $\mathbb{P}(\alpha_m \ge \varepsilon)=1$, for some $\varepsilon \in (0,1)$. Indeed, if this is not the case, then one may introduce another sequence 
\begin{align*}
\widetilde{X}_{t+1,x}^{(m)} = \max (\alpha_{m,t+1},\varepsilon) \widetilde{X}_{t,x}^{(m)} 
 + \beta_m \|\xi_{t+1}\|
\end{align*}
with initial value $\widetilde{X}_{0,x}^{(m)}=\|x\|$ and first hitting time
$\widetilde{\tau} = \min \{t\ge 1: \ \widetilde{X}_{t,x}^{(m)} \le R\}$. Note that $\widetilde{X}_{t,x}^{(m)}\ge \|X_{t,x}^{(m)}\|$ a.s. for all $t\in\N_0$ and, therefore, $\widetilde{\tau}\ge T_x$ a.s. Further, $\varepsilon >0$ may be chosen so small that $\mathbb{E} \log \max (\alpha_m,\varepsilon) <0$.

From now on, assume that $d=1$, $\xi\ge 0$ a.s., $\mathbb{P}(\alpha_m\ge \varepsilon)=1$, $\E\log \alpha_m<0$, and that the initial value $X_{0,x}^{(m)}=x$ is any non-negative number between $0$ and $R$.
Then $X_{t,x}^{(m)}$ is a.s.\ monotone increasing in $x$ for any $t\ge 0$. Then it is enough for \eqref{uniRW} to show that 
\begin{align*}
 {\mathbb E} e^{cT} <\infty
 \end{align*}
for some $c>0$ where $T=T_R = \min \{t>0: \ X_{t,R}^{(m)} \le R\}$.
 
Let $Z_t= \log X_{t,R}^{(m)}$. For any $t\ge 0$, we have
\begin{align*}
Z_{t+1} {\mathbf 1}(T>t)&= \log (\alpha_{m,t} X_{t,R}^{(m)} + \beta_m \xi_{t+1})
{\mathbf 1}(T>t)\\
&\le \log \bigg(\alpha_{m,t}X_{t,R}^{(m)} \bigg(1+\frac{\beta_m\xi_{t+1}}{\varepsilon R}\bigg)\bigg) {\mathbf 1}(T>t)\\
&= 
\bigg(Z_t + \log \alpha_{m,t} + \log \bigg(1+\frac{\beta_m\xi_{t+1}}{\varepsilon R}\bigg)\bigg)
{\mathbf 1}(T>t)\\
&\equiv (Z_t + \psi_{t+1}) {\mathbf 1}(T>t)
\end{align*}
where $\psi_{t+1} = \log \alpha_{m,t} +\log (1+\frac{\beta_m\xi_{t+1}}{\varepsilon R})$, $t\in\N_0$,  are i.i.d.\ random variables. They have a negative mean 
\begin{align*}
{\mathbb E} \psi_{1} \le \mathbb{E} \log \alpha_m + \frac{\beta_m{\mathbb E} \xi_1}{\varepsilon R} <0
\end{align*}
for $R$ large enough since ${\mathbb E} \log\alpha_m<1$. Further, they satisfy ${\mathbb E} e^{\psi_1}= {\mathbb E}\alpha_m (1+ \frac{\beta_m{\mathbb E} \xi_1}{\varepsilon R})<\infty$.
The sequence $\{Z_t\}_{t\in\N_0}$ admits a majorant
$\{\widetilde{Z}_t\}_{t\in\N_0}$ that starts from
$\widetilde{Z}_0=\log R$ and satisfies the recursion
\begin{align*}
\widetilde{Z}_{t+1}	 = \widetilde{Z}_t+\psi_{t+1}, \quad t\in\N_0.
\end{align*}
Then, $\widetilde{T} = \min \{t\ge 1:  \widetilde{Z}_t\le \log R\}\ge T$ a.s.
It is well-known (see e.g.\ \cite{Hey}, Theorem 1 or \cite{Gut}, Theorem 3.2) that the existence of an  exponential moment of
the $\psi$'s and their negative drift imply that $\widetilde{T}$ (and, in turn, $T$) has a finite exponential moment too.

\section{Alternative proof of Theorem \ref{thm:LimitMixing} under stronger moment assumptions by the method of moments}\label{sec:SecondProof}

We provide here an alternative proof of Theorem \ref{thm:LimitMixing}, which contains a number of observations and formulae that are of their own interest, including an inequality for moments.

\subsection{An inequality for moments}

\begin{lemm}\label{lem:InequalityMoments}
Let $Z$ be a random variable such that $\E |Z|^k<\infty$ for all $k\in\N$. If there exists a constant $s\in(0,\infty)$ such that 
\begin{equation}\label{eqn:RecursionMoments}
\E Z^{k} = s (k-1) \E Z^{k-2}
\end{equation}
for all $k\in\N$ with $k\geq 3$, then
\begin{equation}\label{eqn:InequalityMoments}
|\E Z| \leq \frac{\sqrt{2}}{\sqrt{\pi} \sqrt{s}} \E Z^2.
\end{equation}
\end{lemm}

\begin{proof}
We can assume without loss of generality that $\E Z\geq 0$ since we can replace $Z$ by $-Z$, otherwise. For $\ell\in\N$ it follows from \eqref{eqn:RecursionMoments} that
$$
\E Z^{2\ell}= s^{\ell-1} (2\ell-1)\cdot (2\ell-3)\cdot \hdots\cdot 1 \E Z^2= s^{\ell-1} \frac{(2\ell-1)!}{2^{\ell-1}(\ell-1)!} \E Z^2 = s^{\ell-1} \frac{(2\ell-1) (2\ell-2)!}{2^{\ell-1}(\ell-1)!} \E Z^2
$$
and
$$
\E Z^{2\ell-1} = s^{\ell-1} (2\ell-2)\cdot\hdots\cdot 2 \E Z = 2^{\ell-1} s^{\ell-1} (\ell-1)! \E Z.
$$
In the sequel we denote by $\sim$ asymptotic equivalence for $\ell\to\infty$. The Stirling formula leads to
\begin{align*}
\E Z^{2\ell} & \sim  \frac{s^{\ell-1}}{2^{\ell-1}} (2\ell-1) \frac{\sqrt{2\ell-2} (2\ell-2)^{2\ell-2} e^{-(2\ell-2)}}{\sqrt{\ell-1} (\ell-1)^{\ell-1} e^{-(\ell-1)}} \E Z^2\\
& = 2^{\ell-1/2} s^{\ell-1} (2\ell-1) (\ell-1)^{\ell-1} e^{-(\ell-1)} \E Z^2 \\
& \sim 2^{\ell+1/2} s^{\ell-1} (\ell-1)^{\ell} e^{-(\ell-1)} \E Z^2
\end{align*}
and
$$
\E Z^{2\ell-1} \sim 2^{\ell-1} s^{\ell-1} \sqrt{2\pi (\ell-1)} (\ell-1)^{\ell-1} e^{-(\ell-1)} \E Z \sim 2^{\ell-1} s^{\ell-1} \sqrt{2\pi} (\ell-1)^{\ell-1/2} e^{-(\ell-1)} \E Z .
$$
For all $\ell\in\N$ we have
\begin{align*}
\E Z^{2\ell-1} & \leq \E |Z|^{2\ell-1} \leq \E \frac{|Z|^\ell}{(2(\ell-1)s)^{1/4}} (2(\ell-1)s)^{1/4} |Z|^{\ell-1}\\
& \leq \frac{1}{2} \bigg( \frac{\E Z^{2\ell}}{\sqrt{2}\sqrt{\ell-1}\sqrt{s}} + \sqrt{2} \sqrt{\ell-1} \sqrt{s} \E Z^{2\ell-2}\bigg)=:g(\ell).
\end{align*}
The asymptotic formulas formulas from above show that
\begin{align*}
g(\ell) & \sim \frac{1}{2} \bigg( 2^\ell s^{\ell-3/2} (\ell-1)^{\ell-1/2} e^{-(\ell-1)} \E Z^2 + 2^\ell s^{\ell-3/2} \sqrt{\ell-1} (\ell-2)^{\ell-1} e^{-(\ell-2)} \E Z^2 \bigg) \\
& = \frac{1}{2} 2^\ell s^{\ell-3/2}e^{-(\ell-1)} \E Z^2 \bigg(  (\ell-1)^{\ell-1/2} + \sqrt{\ell-1} (\ell-1)^{\ell-1} \bigg(1-\frac{1}{\ell-1} \bigg)^{\ell-1} e \bigg)\\
& \sim 2^\ell s^{\ell-3/2} (\ell-1)^{\ell-1/2} e^{-(\ell-1)} \E Z^2\\
& \sim \frac{\sqrt{2}}{\sqrt{\pi}} \frac{\E Z^2}{\sqrt{s}\E Z} \E Z^{2\ell-1}. 
\end{align*}
Thus, $\E Z^{2\ell-1}\leq g(\ell)$ yields \eqref{eqn:InequalityMoments}.
\end{proof}

\subsection{Second Proof of Theorem \ref{thm:LimitMixing} under stronger moment assumptions}

In this section we prove Theorem \ref{thm:LimitMixing} by the method of moments under the stronger assumptions that
\begin{equation}\label{eqn:ConditionAllMoments}
\E \|\xi\|^k<\infty \quad \text{for all} \quad k\in\N \quad \text{and} \quad \beta_m\sim\frac{1}{\sqrt{m}} \quad \text{as} \quad m\to\infty 
\end{equation}
and that, for all $k\in\N$,
\begin{equation}\label{eqn:ConditionAllMomentsAlpha}
\E \alpha_m^k<1 \quad \text{for all} \quad m\in\N \quad \text{and} \quad 1-\E\alpha_m^k\sim \frac{ka}{m} \quad \text{as} \quad m\to\infty.
\end{equation}

\begin{lemm}\label{prop:Moments}
Assume that \eqref{eqn:ConditionAllMoments} and \eqref{eqn:ConditionAllMomentsAlpha} are satisfied and that $\E \xi=0$. Then, for all $u\in\R^d$, $k\in\N$ and $m\in\N$,
\begin{equation}\label{eqn:RecursiveMoment}
\begin{split}
\E \langle u,Y^{(m)} \rangle^k = & \frac{1}{1-\E \alpha_m^k} \big( -\E\alpha_m^k \E \langle u, Y^{(m)}\rangle^k \mathbf{1}\{Y^{(m)}\in \gamma_m A)\}\\
& \quad \quad \quad \quad + \binom{k}{2} \beta_m^2 \E\langle u,\xi\rangle^2 \E\alpha_m^{k-2} \E \langle u,Y^{(m)}\rangle^{k-2} \mathbf{1}\{Y^{(m)}\notin \gamma_m A\}^{k-2} +R_{k,u}^{(m)}\big)
\end{split}
\end{equation}
with
\begin{equation}\label{eqn:BoundRm}
|R_{k,u}^{(m)}|\leq \beta^3_m\sum_{i=0}^{k-3} \binom{k}{i} \big(1 + |\E \langle u,Y^{(m)}\rangle^i|+ |\E \langle u,Y^{(m)}\rangle^{i+1}|\big) |\E\langle u,\xi\rangle^{k-i}|.
\end{equation}
\end{lemm}

\begin{proof}
It follows from \eqref{eqn:DefinitionYm} and the independence of $\alpha_m$, $Y^{(m)}$ and $\xi$ that
\begin{align*}
\E \langle u,Y^{(m)}\rangle^k & = \E (\alpha_m \langle u,Y^{(m)}\rangle \mathbf{1}\{Y^{(m)}\notin \gamma_m A\} + \beta_m \langle u,\xi\rangle)^k\\
& = \sum_{i=0}^k \binom{k}{i} \E\alpha_m^i \E \langle u,Y^{(m)}\rangle^i \mathbf{1}\{Y^{(m)}\notin \gamma_m A\}^i \beta_m^{k-i} \E \langle u,\xi\rangle^{k-i}.
\end{align*}
Using $\E \xi=0$, we can rewrite this as
\begin{equation}\label{eqn:MomentTemp}
\begin{split}
\E \langle u,Y^{(m)}\rangle^k & = \E\alpha_m^k \E \langle u,Y^{(m)}\rangle^k - \E\alpha_m^k \E \langle u,Y^{(m)}\rangle^k \mathbf{1}\{Y^{(m)}\in \gamma_m A\}\\
& \quad + \binom{k}{2} \beta_m^2 \E \langle u,\xi\rangle^2 \E\alpha_m^{k-2} \E \langle u,Y^{(m)}\rangle^{k-2} \mathbf{1}\{ Y^{(m)}\notin\gamma_m A\}^{k-2} + R_{k,u}^{(m)}
\end{split}
\end{equation}
with
$$
R_{k,u}^{(m)}:= \sum_{i=0}^{k-3} \binom{k}{i} \E\alpha_m^i \E \langle u,Y^{(m)}\rangle^i \mathbf{1}\{Y^{(m)}\notin \gamma_m A\}^i \beta_m^{k-i} \E \langle u,\xi \rangle^{k-i}.
$$
From \eqref{eqn:MomentTemp} we can deduce \eqref{eqn:RecursiveMoment}, while \eqref{eqn:BoundRm} is a consequence of \eqref{eqn:ConditionAllMomentsAlpha}, of $\beta_m\in(0,1)$ and of 
$$
|\E \langle u,Y^{(m)}\rangle^i \mathbf{1}\{Y^{(m)}\notin \gamma_m A\}^i| \leq \E |\langle u,Y^{(m)}\rangle|^i\leq  |\E \langle u,Y^{(m)}\rangle^i|+1 + |\E \langle u,Y^{(m)}\rangle^{i+1}|
$$
for $i\in\{1,\hdots,k-3\}$, where the last inequality follows from considering even and odd $i$ separately.
\end{proof}

\begin{prop}\label{prop:AsymptoticsMoments2}
Let $u\in\R^d$ and assume that \eqref{eqn:ConditionAllMoments} and \eqref{eqn:ConditionAllMomentsAlpha} are satisfied and that $\E \xi=0$. If
$$
\E\langle u, Y^{(m)} \rangle \to \mu_{1,u}\in\R \quad \text{and} \quad \E \langle u,Y^{(m)}\rangle^2\to\mu_{2,u}\in[0,\infty) \quad \text{as} \quad m\to\infty,
$$
then all limits 
$$
\mu_{k,u}:=\lim_{m\to\infty} \E \langle u,Y^{(m)}\rangle^k
$$
with $k\in\N$ and $k\geq 3$ exist as finite numbers and
\begin{equation}\label{eqn:RecursionMomentsYm}
\mu_{k,u}=\frac{k-1}{2a} \E\langle u,\xi\rangle^2 \mu_{k-2,u}
\end{equation}
for $k\in\N$ with $k\geq 3$.
\end{prop}

\begin{proof}
We prove the statement for $k\in\N$ with $k\geq 3$ under the assumption that the assertion is true for all moments of order smaller than $k$. Obviously this assumption is satisfied for $k=3$ and we can iterate the argument. From Lemma \ref{prop:Moments} it follows that
\begin{equation}\label{eqn:DecompositionMoment}
\begin{split}
\E \langle u,Y^{(m)} \rangle^k = & \frac{1}{1-\E\alpha_m^k} \bigg( -\E\alpha_m^k \E \langle u,Y^{(m)}\rangle^k \mathbf{1}\{Y^{(m)}\in\gamma_m A\}\\
& \quad \quad \quad \quad + \binom{k}{2} \beta_m^2 \E\langle u,\xi\rangle^2 \E\alpha_m^{k-2} \E \langle u,Y^{(m)}\rangle^{k-2} \mathbf{1}\{Y^{(m)}\notin \gamma_m A\} +R_{k,u}^{(m)}\bigg).
\end{split}
\end{equation}
Because of \eqref{eqn:ConditionAllMomentsAlpha}, the prefactor behaves as $m/(ka)$ as $m\to\infty$. Lemma \ref{lem:BoundMomentInterval} leads to
$$
m \E\alpha_m^k |\E \langle u,Y^{(m)}\rangle^k \mathbf{1}\{Y^{(m)}\in\gamma_m A\}| \leq m \E\alpha_m^k \|u\|^k \bigg( 1 + \frac{\overline{r}^k}{\underline{r}^k} \bigg)  \beta_m^k \frac{\E\|\xi\|^k}{\E \tau^{(m)}},
$$
whence the first summand in \eqref{eqn:DecompositionMoment} vanishes as $m\to\infty$. Here, we used that $\E \tau^{(m)}\geq 1$. From the assumed convergence of all moments of order less than $k$ and \eqref{eqn:BoundRm} we obtain that
$$
\limsup_{m\to\infty} m|R_{k,u}^{(m)}| \leq \lim_{m\to\infty} m \beta_m^3 \sum_{i=0}^{k-3} \binom{k}{i} \big( 1 + |\E \langle u,Y^{(m)} \rangle^i| + |\E \langle u,Y^{(m)}\rangle^{i+1}|\big) |\E\langle u,\xi \rangle^{k-i}| =0
$$
so that the last term in \eqref{eqn:DecompositionMoment} goes to zero as $m\to\infty$. Lemma \ref{prop:Moments} yields that
$$
\E \langle u,Y^{(m)} \rangle = -\frac{\E\alpha_m}{1-\E\alpha_m} \E \langle u,Y^{(m)} \rangle \mathbf{1}\{Y^{(m)}\in \gamma_m A\}.
$$
Together with $\E \langle u,Y^{(m)} \rangle \to\mu_{1,u}$ and $\E\alpha_m/(1-\E\alpha_m)\to\infty$ as $m\to\infty$, we deduce that
$$
\lim_{m\to\infty} \E \langle u,Y^{(m)} \rangle \mathbf{1}\{Y^{(m)}\in \gamma_m A\}=0. 
$$
Thus, for $k=3$ we have
\begin{align*}
& \lim_{m\to\infty} m \binom{k}{2} \beta_m^2 \E\langle u,\xi\rangle^2 \E\alpha_m^{k-2} \E \langle u,Y^{(m)}\rangle^{k-2} \mathbf{1}\{Y^{(m)}\notin \gamma_m A\} \\
& = 3 \E\langle u,\xi\rangle^2 \lim_{m\to\infty} \E \langle u,Y^{(m)}\rangle - \E \langle u,Y^{(m)}\rangle \mathbf{1}\{Y^{(m)}\in \gamma_m A\} \\
& = 3 \E\langle u,\xi\rangle^2 \mu_{1,u} = \binom{k}{2} \E\langle u,\xi\rangle^2 \mu_{k-2,u}.
\end{align*}
For $k>3$, it follows from Lemma \ref{lem:BoundMomentInterval} that
\begin{align*}
& \limsup_{m\to\infty} m \binom{k}{2} \beta_m^2 \E\langle u,\xi \rangle^2 \E\alpha_m^{k-2} \big|\E \langle u,Y^{(m)}\rangle^{k-2} \mathbf{1}\{Y^{(m)}\in \gamma_m A\}\big|\\
& \leq \binom{k}{2} \E\langle u,\xi\rangle^2 \limsup_{m\to\infty} \E |\langle u,Y^{(m)}\rangle|^{k-2} \mathbf{1}\{Y^{(m)}\in \gamma_m A\} \\
& \leq \binom{k}{2} \E\langle u,\xi\rangle^2 \|u\|^{k-2} \bigg( 1+ \frac{\overline{r}^{k-2}}{\underline{r}^{k-2}} \bigg) \frac{\E\|\xi\|^{k-2}}{\E \tau^{(m)}} \lim_{m\to\infty} \beta_m^{k-2}=0,
\end{align*}
whence
\begin{align*}
& \lim_{m\to\infty} m\binom{k}{2} \beta_m^2 \E\langle u,\xi\rangle^2 \E\alpha_m^{k-2} \E \langle u,Y^{(m)}\rangle^{k-2} \mathbf{1}\{Y^{(m)}\notin \gamma_m A\}\\
& = \lim_{m\to\infty} \binom{k}{2} \E\langle u,\xi\rangle^2 \E \langle u,Y^{(m)}\rangle^{k-2} = \binom{k}{2} \E\langle u,\xi\rangle^2 \mu_{k-2,u}.
\end{align*}
Combining these computations, we obtain
$$
\mu_{k,u} = \lim_{m\to\infty} \E \langle u,Y^{(m)}\rangle^k = \frac{1}{ka} \binom{k}{2} \E\langle u,\xi\rangle^2 \mu_{k-2,u} = \frac{k-1}{2a} \E\langle u,\xi\rangle^2 \mu_{k-2,u},
$$
which completes the proof.
\end{proof}

\begin{proof}[Proof of Theorem \ref{thm:LimitMixing} under assumptions \eqref{eqn:ConditionAllMoments} and \eqref{eqn:ConditionAllMomentsAlpha}]

As in the proof by characteristic functions we obtain that $\{Y^{(m)}\}_{m\in\N}$ has a subsequence that converges in distribution. We assume for simplicity that this subsequence is the sequence itself and denote the limiting random vector by $Y$. Let $v\in\R^d$ be fixed. By \eqref{EY} and Proposition \ref{prop:2ndMoment} we obtain
$$
\mu_{1,v}:=\lim_{m\to\infty} \E \langle v,Y^{(m)}\rangle= \langle v,\mu\rangle \quad \text{and} \quad \mu_{2,v}:=\lim_{m\to\infty} \E \langle v,Y^{(m)}\rangle^2= p \frac{\E\langle v,\xi\rangle^2}{2a}.
$$
It follows from Proposition \ref{prop:AsymptoticsMoments2} that the limits
$$
\mu_{k,v}:=\lim_{m\to\infty} \E\langle u, Y^{(m)} \rangle^k
$$
for $k\in\N$ with $k\geq 3$ exist and satisfy \eqref{eqn:RecursionMomentsYm}. Since the moments of $\{\langle v,Y^{(m)}\rangle \}_{m\in\N}$ converge to those of $\langle v,Y \rangle$, we have $\E \langle v,Y\rangle^k=\mu_{k,v}$ for $k\in\N$. We can apply Lemma \ref{lem:InequalityMoments} to deduce that
$$
|\langle v,\mu \rangle| = |\mu_{1,v}| \leq \frac{\sqrt{2}}{\sqrt{\pi}} \frac{\sqrt{2a}}{\sqrt{\E\langle v,\xi\rangle^2}} \mu_{2,v} = \frac{p}{\sqrt{\pi a}} \sqrt{\E\langle v,\xi\rangle^2} = \frac{p}{\sqrt{\pi a}} \sqrt{v^T\Sigma v}.
$$
This shows that the random variable $B_{2,v}$ with the probabilities given in \eqref{defB} is well defined.

We define 
$$
Y_v = \frac{\sqrt{v^T\Sigma v}}{\sqrt{2a}} B_1 B_{2,v} |N| = \frac{\sqrt{\E \langle v, \xi \rangle^2}}{\sqrt{2a}} B_1 B_{2,v} |N|
$$
with $B_1$, $B_{2,v}$ and $N$ as in Theorem \ref{thm:LimitMixing} so that
$$
\E Y_v =\langle v,\mu\rangle = \mu_{1,v} \quad \text{and} \quad \E Y_v^2 = p\frac{\E\langle v,\xi\rangle^2}{2a}=\mu_{2,v}. 
$$
Because of
$$
\E |N|^k = (k-1) \E |N|^{k-2}
$$
for $k\in\N$ with $k\geq 3$, we have that
$$
\E Y_v^k = (k-1) \frac{\E \langle v,\xi \rangle^2}{2a} \E Y_v^{k-2}
$$
for $k\in\N$ with $k\geq 3$. Since this is the same recursive formula as \eqref{eqn:RecursionMomentsYm} in Proposition \ref{prop:AsymptoticsMoments2}, we obtain that $\E Y_v^k=\mu_{k,v}$ for $k\in\N$. The naive bound
$$
|\mu_{k,v}| \leq k!\max\bigg\{1, \frac{\E\langle v, \xi \rangle^2}{2a}\bigg\}^k \max\{|\mu_{1,v}|,|\mu_{2,v}|\}
$$
for $k\in\N$ implies 
$$
\limsup_{k\to\infty} \frac{|\mu_{k,v}|^{1/k}}{k}<\infty,
$$
whence the distributions of $Y_v$ and $\langle v,Y\rangle$ are completely determined by their moments (see \cite[Chapter 2, § 12, Theorem 7]{Shiryaev}). Thus, the method of moments yields the distributional convergence
$$
\langle v, Y^{(m)} \rangle\overset{d}{\longrightarrow} \frac{\sqrt{v^T\Sigma v}}{\sqrt{2a}} B_1 B_{2,v} |N| \quad \text{as} \quad m\to\infty.
$$ 
Together with the Cramer-Wold device we obtain \eqref{eqn:LimitMixing} and \eqref{eqn:vProjection}. From the recursive formula \eqref{eqn:RecursionMomentsYm} for the moments of $\langle v, Y\rangle$, one can directly compute the characteristic function of $Y$ and, thus, $\varphi_Z$ in \eqref{eqn:Formula_phi}. From the characteristic function one can derive the density in \eqref{densZ} as in the original proof.
\end{proof}

\section*{Acknowledgements}

The authors are thankful to Gavin Gibson for introducing us to an interesting model,
and to Onno Boxma, Charles Goldie, Alexander Lindner, Thomas Mikosch, Vitali Wachtel and Stan Zachary for helpful comments.


\begin{thebibliography}{}

\bibitem{AlsIksRoe} G.~Alsmeyer, A.~Iksanov and U.~R\"{o}ssler (2009), On distributional
properties of perpetuities, {\em Journal of Theoretical Probability}, {\bf 22}, 666--682.

\bibitem{Asm} S.~Asmussen (2003), {\em Applied Probability and Queues}, 2nd Edition, Springer. 

\bibitem{Bor0} A.A. Borovkov (1976), {\em Stochastic Processes in Queueing Theory}, Wiley.

\bibitem{Bor} A.~Borovkov and S.~Foss (1992), Stochastically
recursive sequences and their generalisations, {\em Siberian Advances in Mathematics}, {\bf 2}, 16--81.

\bibitem{BMR} O.~Boxma, M.~Mandjes and J.~Reed (2016),
On a Class of reflected $AR(1)$ processes, {\em Journal of Applied Probability}, {\bf 53}, 818--832.

\bibitem{Bra} A.~Brandt (1986), 
The stochastic equation $Y_{n+1}=A_nY_n+B_n$ with stationary coefficients, {\em Advances in
Applied Probability}, {\bf 18}, 211--220.

\bibitem{Bur3} 
D.~Buraczewski, J.F.~Collamore, E.~Damek and J.~Zienkiewicz (2016),
Large deviation estimates for exceedance times of perpetuity sequences and their dual processes, {\em The Annals of Probability}, {\bf 44}, 3688--3739.

\bibitem{Bur1} D.~Buraczewski, E.~Damek and T.~Mikosch (2016), {\em Stochastic models with power-law tails. The
equation $X=AX+B$},  Springer.

\bibitem{Bur2} D.~Buraczewski and A.~Iksanov (2015), Functional limit theorems for divergent perpetuities in the
contractive case, {\em Electronic Communications in 
Probability}, {\bf 20}.


\bibitem{DaMi} R.A.~Davis and T.~Mikosch (1998),
Gaussian likelihood-based inference for non-invertible
$MA(1)$ processes with $S\alpha S$ noise,
{\em Stochastic Processes and their Applications}, {\bf 77}, 99--122.

\bibitem{Emb} P.~Embrechts and C.~Goldie (1994), Perpetuities and Random Equations, in {\em Asymptotic Statistics}, Springer, 
75--86.

\bibitem{Erh} T.~Erhardsson (2014), Conditions for convergence of random coefficient $AR(1)$ processes and
perpetuities in higher dimensions, {\em Bernoulli}, {\bf 20},
990--1005.

\bibitem{Gut} A.~Gut (2009), {\em Stopped Random Walks}, 2nd Edition, Springer.

\bibitem{Hey} C.C.~Heyde (1964), Two probability theorems and their applications to some
first passage problems, {\em Journal of the Australian Mathematical Society}, {\bf 4}, 214--222.

\bibitem{Kes} H.~Kesten (1973), Random difference equations and renewal theory for products of random matrices, 
{\em Acta Mathematica}, {\bf 131}, 207--248.

\bibitem{Kif} Y.~Kifer (1986), {\em Ergodic Theory of Random Transformations}, Birkh\"{a}user.

\bibitem{Kol} B.~Kolodziejek (2018), On perpetuites with light tails, {\em Advances in Applied Probability}, {\bf 50}, 1119--1154. 

\bibitem{Rah} T.~Lange and A.~Rahbek, An Introduction to Regime Switching Time Series Models, 
 in: {\em Handbook on Financial Time Series},
T.G.~Andersen, R.A.~Davis, J.-P.~Kreiss, T.~Mikosch (Editors),
Springer, 2009,  871--888.

\bibitem{Neveu} J.~Neveu (1975), {\em Discrete-parameter martingales}, Elsevier.

\bibitem{Shiryaev} A.N.~Shiryaev (1996), {\em Probability}, 2nd Edition, Springer.

\bibitem{Vervaat1979} W.~Vervaat (1979), On a stochastic difference equation and a representation of nonnegative infinitely divisible random variables, {\em Advances in Applied Probability}, {\bf 11}, 750--783.

\bibitem{Zer} M.~Zerner (2018), Recurrence and transience of contractive autoregressive processes and related Markov chains,
{\em Electronic Journal of Probability}, {\bf 23}, paper no. 27, 24pp.

\end{thebibliography}
\end{document}